\documentclass[a4paper,11pt]{article}
\usepackage{fullpage,amssymb,latexsym,amsmath,amsthm,}
\usepackage[latin1]{inputenc}

\newcommand{\dN}{\ensuremath{\mathbb{N}}} 
\newcommand{\dR}{\ensuremath{\mathbb{R}}} 

\newcommand{\ind}{\mathrm{1}\hskip -3.2pt \mathrm{I}} 
\newcommand{\ent}{\mathbf{Ent}} 
\newcommand{\var}{\mathbf{Var}} 

\newtheorem{theorem}{Theorem}
\newtheorem{proposition}[theorem]{Proposition}
\newtheorem{lemma}[theorem]{Lemma}
\newtheorem{corollary}[theorem]{Corollary}

\theoremstyle{definition} \newtheorem{definition}{Definition}

\theoremstyle{remark}
\newtheorem{remark}[theorem]{Remark}
\newtheorem{example}[theorem]{Example}

\begin{document}   

 \title{Modified logarithmic Sobolev inequalities on $\mathbb R$}
 \author{F. Barthe and C. Roberto}
 \maketitle
 \date{}

\begin{abstract}
  We provide a sufficient condition for a measure on the real line to
  satisfy a modified logarithmic Sobolev inequality, thus extending the criterion of Bobkov
  and Götze. Under mild assumptions
  the condition is also necessary. Concentration inequalities are derived.
   This  completes the picture given in recent contributions by Gentil, Guillin and Miclo.
\end{abstract}

\section{Introduction}
 In this paper we are interested in Sobolev type inequalities satisfied by probability
measures. It is well known that they allow to describe their concentration properties
as well as the regularizing effects of associated semigroups. Several books are available on these
topics and we refer to them for more details (see e.g. \cite{Ane,ledoCMLS}).
Establishing such inequalities is a difficult task in general, especially in high dimensions.
However, it is very natural to investigate such inequalities for measures on the real line. Indeed many
high dimensional results are obtained by induction on dimension, and having a good knowledge
of one dimensional measures becomes crucial. Thanks to Hardy-type inequalities, it is possible
to describe very precisely the measures on the real line which satisfy certain Sobolev inequalities.
Our goal here is to extend this approach to the so-called modified logarithmic Sobolev inequalities.
They are introduced below.

\medskip

Let $\gamma$ denote the standard Gaussian probability measure on $\mathbb R$.
 The Gaussian logarithmic Sobolev asserts that for every smooth $f:\mathbb R\to \mathbb R$
 $$ \ent_\gamma(f^2) \le 2 \int (f')^2 d\gamma,$$
where the entropy functional with respect to a probability measure $\mu$ is defined by
$$\ent_\mu(f)=\int f\log f \, d\mu
-\left(\int f\, d\mu\right)  \log\left(\int f \, d\mu \right).$$
This famous inequality implies the Gaussian concentration inequality, as well as hypercontractivity
and entropy decay along the Ornstein-Uhlenbeck semigroup. Since the logarithmic-Sobolev inequality
implies a sub-Gaussian behavior of tails, it is not verified for many measures and one has to
consider weaker Sobolev inequalities. In the case of the symmetric exponential measure $d\nu(t) =e^{-|t|}dt/2$,
an even more classical fact is available, namely a Poincaré or spectral gap inequality: for every
smooth function $f$:
\begin{equation}\label{eq:Pe}
 \var_\nu (f)\le 4 \int (f')^2 d\nu.
\end{equation}
This property implies an exponential concentration inequality, as noted by Gromov and Milman \cite{gromm83taii}, as well
as a fast decay of the variance along the corresponding semigroup. If one compares to the log-Sobolev
inequality, the spectral gap inequality differs by its left side only. In order to describe more
precisely the concentration phenomenon  for product of exponential measures, and to recover a celebrated
result by Talagrand \cite{tala91niic},
 Bobkov and Ledoux \cite{bobkl97pitc} introduced a so-called modified logarithmic Sobolev inequality
for the exponential measure. Here the entropy term remains but the term involving the derivatives is changed.
Their result asserts that every smooth $f:\mathbb R\to \mathbb R$ with $ |f'/f|\le c<1$ verifies
\begin{equation}\label{eq:LSe}
 \ent_\nu(f^2) \le \frac{2}{1-c} \int (f')^2 d\nu.
\end{equation}
The latter may be rewritten as
\begin{equation}\label{eq:LMSe}
 \ent_\nu(f^2) \le \int H\left(\frac{f'}{f}\right)f^2 d\nu,
\end{equation}
where  $H(t)=2t^2/(1-c)$ if $|t|\le c$ and $H(t)=+\infty$ otherwise.
Such general modified log-Sobolev inequalities have been established
by Bobkov and Ledoux \cite{bobkl00bmbl} for the probability measures $d\nu_p(t)=e^{-|t|^p}dt/Z_p$, $t\in \mathbb R$
in the case $p>2$ (a more general result is valid for measures $e^{-V(x)}dx $ on $\mathbb R^n$
where $V$ is strictly uniformly convex). These measures satisfy a modified log-Sobolev inequality
with function $H(t)=c_p |t|^{q}$ where $q=p/(p-1) \in [1,2]$ is the dual exponent of $p\ge 2$.
The inequality can be reformulated as $\ent_{\nu_p}(|g|^q) \le \tilde{c}_p \int |g'|^q d\nu_p$.
These $q$-log-Sobolev inequalities are studied in details by Bobkov and Zegarlinski in \cite{bobkz05ebi}.

 The case $p\in(1,2)$ is more delicate: the inequality cannot hold with  $H(t)=c_p |t|^{q}$
since this function is too small close to zero. Indeed for $f=1+\varepsilon g$ when  $g$ is bounded
and $\varepsilon$ very small, the left hand side of \eqref{eq:LMSe} is equivalent to $\varepsilon^2 \var_\nu(g)$
whereas the right hand side is comparable to $\int H(\varepsilon g') \,d\nu$. Hence $H(t)$ cannot be much smaller
than $t^2$ when $t$ goes to zero. If it compares to $t^2$ then in the limit one recovers a spectral gap inequality.
 Gentil, Guillin and Miclo \cite{gentilgm} established a modified log-Sobolev inequality
for $\nu_p$ when $p\in (1,2)$, with a function $H_p(t)$ comparable to $k_p \max(t^2, |t|^q)$. In the subsequent
paper \cite{gentilgm2} they extend their method to even log-concave measures on the line, with tail behavior between
exponential and Gaussian. Their method is rather involved. It relies on classical Hardy types
inequalities, adapted to inequalities involving terms as $\int (f')^2 d\mu$, where $\mu$ is carefully chosen.

Our alternative approach  is to develop Hardy type methods directly for inequalities involving  terms as $\int H(f'/f) f^2 d\mu$.
This is done abstractly in Section~2, but more work is needed to present the results in an explicit and workable form.
Section~3 provides a simple sufficient condition for a measure to satisfy a modified log-Sobolev inequality
with function $H(t)=k_p \max(t^2, |t|^q)$ for $p\in (1,2)$, and recovers in a soft way the result of  \cite{gentilgm}.
Under mild assumptions, the condition is also necessary and we have a reasonable estimate of the
best constant in the inequality. Next in Section~4 we consider the same problem for general convex functions $H$. The
approach remains rather simple, but technicalities are more involved. However Theorem~\ref{th:critphi}
provides a neat sufficient condition, which recovers the result of  \cite{gentilgm2} for log-concave measures but also
applies without this restriction. Under a few more assumptions, our sufficient condition is also necessary.
In Section~5 we describe concentration consequences of modified logarithmic Sobolev inequalities,
obtained by the Herbst method.

Logarithmic Sobolev inequalities are known to imply inequalities between transportation cost and
entropy \cite{ottov00gitl,bobkgl01hhje}. Our criterion can be compared with the one recently derived by
Gozlan \cite{gozl06ctti}. It confirms that modified logarithmic Sobolev inequalities are strictly stronger
than the corresponding transportation cost inequalities, as
discovered by Cattiaux and Guillin \cite{cattg06qtci} for the classical logarithmic Sobolev inequality and Talagrand's
transportation cost inequality. For log-concave measures on $\mathbb R$ the results of Gozlan yield precise
modified logarithmic Sobolev inequalities. By different methods, based on isoperimetric inequalities, Kolesnikov \cite{kole06mlsi}
recently established more general modified $F$-Sobolev inequalities for log-concave probability measures on $\mathbb R^n$.

\bigskip
We end this introduction by setting the notation.
It will be convenient to work with locally Lipschitz functions $f:\mathbb R^d \to \mathbb R$,
for which the norm of the gradient (absolute value of the derivative when $d=1$) can be defined as a
whole by
$$ |\nabla f|(x)=\lim_{r\to 0^+}\sup_{y; \; |x-y|\le r}\frac{|f(x)-f(y)|}{|x-y|},$$
where the denominator is the Euclidean norm of $x-y$. By Rademacher's theorem, $f$ is Lebesgue almost
everywhere differentiable, and at these points the above notion coincides with the Euclidean norm of the
gradient of $f$.

We recall that a Young function is an even convex function $\Phi:\mathbb R\to [0,+\infty)$ with
$\Phi(0)=0$ and $\lim_{x\to +\infty }\Phi(x)=+\infty $. Following \cite{rao-ren} we say that $\Phi$
is a nice Young function if it also verifies $\Phi'(0)=0$, $\lim_{x\to+\infty }\frac{\Phi(x)}{x}=+\infty$
and   vanishes only at $0$. We refer to the Appendix for more details about these functions and
their Legendre transforms.

Given a nice Young function $\Phi : \dR^+ \rightarrow \dR^+$ we define its  \emph{modification}
\begin{equation} \label{eq:m}
  \begin{array}{rcl}
    H_\Phi : \dR & \rightarrow & \dR^+ \\
    x & \mapsto & x^2 \ind_{[0,1]} + \frac{\Phi(|x|)}{\Phi(1)}\ind_{]1,\infty)} .
  \end{array}
\end{equation}
%
A probability measure $\mu$ on $\dR$ satisfies a
modified logarithmic Sobolev inequality with function
$H_\Phi$, if there exists some constant $\kappa \in (0,\infty)$
such that every locally Lipschitz $f : \dR \rightarrow \dR$ satisfies
$$
\ent_\mu (f^2) \leq \kappa \int H_\Phi \left( \frac{f'}{f} \right) f^2 d\mu .
$$
We consider functions $\Phi$ such that $\Phi(x)\ge cx^2$ for $x\ge 1$, hence the inequalities
we study are always weaker than the classical logarithmic Sobolev inequality.
On the other hand, as recalled in the introduction,  they  imply the Poincaré Inequality.




\section{Hardy inequalities on the line}
In this section we show how the modified log-Sobolev inequality can be
addressed by Hardy type inequalities.  We refer to the book \cite{Ane}
for the history of the topic.  The extension of Hardy's inequalities
to general measures, due to Muckenhoupt \cite{muck72hiw}, allowed
recent progress in the understanding of several functional
inequalities on the real line.
We recall it below:
\begin{theorem}\label{th:hardy}
   Let $\mu,\nu$ be Borel measures on $\mathbb R^+$ and $p>1$. Then
   the best constant $A$ such that every locally Lipschitz function $f$ verifies
   $$
   \int_{[0,+\infty)} |f-f(0)|^p d\mu \le A \int_{[0,+\infty)}
   |f'|^p d\nu $$
   is finite if and only if
   $$
   B:=\sup_{x>0} \mu\big([x,+\infty) \big) \left( \int_0^x
      \frac{1}{n^{\frac{1}{p-1}}} d\nu\right)^{p-1}$$
   is finite. Here
   $n$ is the density of the absolutely continuous part of $\nu$.
   Moreover, when it is finite $B\le A\le \frac{p^p}{(p-1)^{p-1}}B$.
\end{theorem}
As an easy consequence, one gets a characterization of
measures satisfying a spectral gap inequality together with a good
estimate of the optimal constant  (see e.g \cite{Ane}). The next statement also
gives an improved lower bound on the best constant $C_P$ recently obtained by Miclo
\cite{micl06qbhp}.
\begin{theorem}\label{th:hpoinc}
   Let $\mu$ be a probability measure on $\mathbb R$ with median $m$
   and let $d\nu(t)=n(t) \, dt$ be a measure on $\mathbb R$.  The best
   constant $C_P$ such that every locally Lipschitz $f:\mathbb R\to \mathbb R$
   verifies
   \begin{equation}\label{eq:poincare}
   \var_\mu(f)\le C_P \int (f')^2 d\nu
   \end{equation}
   verifies $ \max(B_+,B_-) \le C_P \le 4 \max(B_+,B_-)$, where
   $$
   B_+=\sup_{x>m} \mu\big([x,+\infty) \big) \int_m^x \frac{1}{n},
   \quad B_-=\sup_{x<m} \mu\big((-\infty,x] \big) \int_x^m \frac{1}{n}
   \cdot$$
\end{theorem}
Bobkov and Götze \cite{bobkg99eitc} used Hardy inequalities to obtain
a similar result for the best constant in logarithmic Sobolev
inequalities: they showed that up to numerical constants, the best
$C_{LS}$ such that for all locally Lipschitz $f$
$$
\ent_\mu(f)\le C_{LS} \int (f')^2 d\nu,$$
is the maximum of
$$
\sup_{x>m} \mu\big([x,+\infty)
\big)\log\left(\frac{1}{\mu\big([x,+\infty) \big)}\right) \int_m^x
\frac{1}{n}$$
and of the corresponding term involving the left side of
the median. In \cite{bartr03sipm}, we improved their method and
extended it to inequalities interpolating between Poincaré and
log-Sobolev inequalities (but involving $\int (f')^2d\nu$).

Using classical arguments (see e.g. the Appendix of
\cite{micl06qbhp}) it is easy to see that the Poincaré, the
logarithmic Sobolev and the modified logarithmic Sobolev constants
are left unchanged if one restrict oneself to the absolutely
continuous part of the measure $\nu$ in the right hand side. So,
without loss of generality, in the sequel we will always assume that
$\nu$ is absolutely continuous with respect to the Lebesgue measure.

The next two statements show that similar results hold for modified
log-Sobolev inequalities provided one replaces the term $\int_m^x 1/n$
by suitable quantities. Obtaining workable expressions for them is not
so easy, and will be addressed in the next sections.

\begin{proposition}\label{prop:maj}
   Let $\mu$ be a probability measure with median $m$ and $\nu$ a
   non-negative measure, on $\dR$.
  Assume that $\nu$ is absolutely continuous with respect to Lebesgue measure and
   that the following Poincaré inequality is satisfied: for all locally Lipschitz $f$
  $$ \var_\mu(f) \le C_P \int (f')^2 d\nu.$$

 Let $\Phi$ be a nice Young
   function such that $\Phi(t)/t^2$ is non-decreasing for $t>0$.
   Define for $x >m$ the number $\alpha_x^+$ and for $x<m$ the number
   $\alpha_x^-$ as follows
   $$
   \alpha_x^+:= \inf \left\{ \int_m^x \Phi \left(\frac{f'}{f}
      \right) f^2 d\nu , f \mbox{non-decreasing}, f(m)=1, f(x)=2
   \right\},
   $$
   $$
   \alpha_x^-:= \inf \left\{ \int_x^m \Phi \left(\frac{f'}{f}
      \right) f^2 d\nu , f \mbox{non-increasing}, f(x)=2, f(m)=1
   \right\} .
   $$
   Denote
\begin{eqnarray*}
B^+ (\Phi): & = &
\sup_{x >m} \mu([x,\infty)) \log \left(\frac{1}{\mu([x,\infty))}\right) \frac{1}{\alpha_x^+} , \\
B^- (\Phi): & = &
\sup_{x <m} \mu((-\infty,x]) \log \left(\frac{1}{\mu((-\infty,x])}\right) \frac{1}{\alpha_x^-} .
\end{eqnarray*}
Then for any for any locally Lipschitz $f: \dR \rightarrow \dR$
\begin{equation*}
\ent_\mu (f^2) \leq \Big(235 C_P+8\Phi(1)\max\big(B_+(\Phi),B_-(\Phi)\big)\Big)
 \int H_\Phi \left(\frac{f'}{f} \right) f^2 d\nu .
\end{equation*}
\end{proposition}

\begin{proof}
   In the above statement, there is nothing canonical about 2 in the
   definition of $\alpha^+_x$ and $\alpha_x^-$. We could  replace
   it by a parameter $\sqrt{\rho}>1$. Optimising over $\rho$ would
   yield non-essential improvements in the results of this paper.
   However, for this proof we keep the parameter, as we find it
   clearer like this. We set $\rho=4$ and any value stricty bigger
   than 1 would do.

   Without loss of generality we start with a non-negative function
   $f$ on $\mathbb R$.  We consider the associated function
   \begin{eqnarray*}
      g(x)&= &f(m)+\int_{m}^x f'(u) \ind_{f'(u)>0}\, du \quad \mathrm{if} \quad x\ge m \\
        g(x)&= &f(m)+\int_{m}^x f'(u) \ind_{f'(u)<0}\, du \quad \mathrm{if} \quad x< m.
   \end{eqnarray*}
   We follow the method of Miclo-Roberto \cite[Chapter 3]{robe} (see also Section~5.5 of
   \cite{bartcr04iibe} where it is extended). We will omit a few details, which are available
   in these references. We introduce  for $x,t >0$, $\Psi_t(x)=x \log(x/t)-(x-t)$.
By convexity of the function $x\log x$ it is easy to check that
$$ \ent_\mu(f^2)  = \int \Psi_{\mu(f^2)} (f^2) d\mu =\inf_t  \int \Psi_{t} (f^2) d\mu
    \le  \int \Psi_{\mu(g^2)} (f^2) d\mu.$$
Defining $\Omega:=\{x; \; f^2(x)\ge 2\rho\, \mu(g^2) \}$, we get
\begin{equation}\label{eq:start}
 \ent_\mu(f^2)    \le  \int_{\Omega^c}  \Psi_{\mu(g^2)} (f^2) d\mu
  + \int_{\Omega\cap [m,+\infty)}  \Psi_{\mu(g^2)} (f^2) d\mu
 + \int_{\Omega\cap (-\infty,m]}  \Psi_{\mu(g^2)} (f^2) d\mu .
\end{equation}
The first term is bounded as follows.
One can check that for any $x \in[0,\sqrt{2 \rho}t]$, it holds $\Psi_{t^2}(x^2) \leq (1+\sqrt{2
  \rho})^2(x-t)^2$. Thus
\begin{eqnarray*}
&&\int_{\Omega^c} \Psi_{\mu(g^2)} (f^2) d\mu
 \leq
(1+\sqrt{2 \rho})^2 \int_{\Omega^c} \left( f-\sqrt{\mu(g^2)} \right)^2 d\mu \\
& \le  &
2(1+\sqrt{2 \rho})^2 \int  \big( f-g\big)^2 d\mu + 2(1+\sqrt{2 \rho})^2
\int  \big(g - \sqrt{\mu(g^2)} \big)^2 d\mu
\end{eqnarray*}
The last term of the above expression is bounded from above by applying the Poincaré inequality to $g$.
Using the definition of $g$ and applying Hardy's inequality on $(-\infty,m]$ and $[m,+\infty)$ allows
to upper bound the term $\int (f-g)^2 d\mu$. By Theorems~\ref{th:hardy} and \ref{th:hpoinc}
 the best constants in Hardy inequality
compare to the Poincaré constant. Finally one gets
$$   \int_{\Omega^c}  \Psi_{\mu(g^2)} (f^2) d\mu \le 16(1+\sqrt{2\rho})^2 C_P \int (f')^2 d\nu.$$
The second term in \eqref{eq:start} is
$$ \int_{[m,+\infty)\cap \{f^2\ge 2\rho\,\mu(g^2)\}}
 \left(f^2 \log\Big(\frac{f^2}{\mu(g^2)}\Big)-(f^2-\mu(g^2))\right)d\mu$$
$$ \le  \int_{[m,+\infty)\cap \{g^2\ge 2\rho\,\mu(g^2)\}} g^2 \log\Big(\frac{g^2}{\mu(g^2)}\Big)d\mu=
 \int_{\Omega_1} g^2 \log\Big(\frac{g^2}{\mu(g^2)}\Big)d\mu,$$
where we have set for  $k \in \dN$, $\Omega_k:=\left\{x\ge m;\; g^2(x) \geq 2 \rho^k \mu(g^2)
\right\}$.  Since $g$ is non-decreasing on the right of $m$,
 we have $\Omega_{k+1} \subset
\Omega_k=[a_k,\infty)$ for some $a_k \geq m$. Also by Markov's inequality $\mu(\Omega_k)\le 1/(2\rho^k)$.
Furthermore, on $\Omega_k \setminus \Omega_{k+1}$, $2 \rho^k \mu(g^2)
\leq g^2 <2 \rho^{k+1} \mu(g^2)$. Thus we have
\begin{eqnarray*}
\int_{\Omega_1} g^2 \log \frac{g^2}{\mu(g^2)} d\mu
& = &
\sum_{k \geq 1} \int_{\Omega_k \setminus \Omega_{k+1}}
g^2 \log \frac{g^2}{\mu(g^2)} d\mu \\
& \leq &
\sum_{k \geq 1} \mu(\Omega_k) 2 \rho^{k+1} \mu(g^2) \log (2 \rho^{k+1}) \\
& \leq &
2\sum_{k \geq 1} \mu(\Omega_k) 2 \rho^{k+1} \mu(g^2) \log (2 \rho^{k})\\
& \leq &
2 \sum_{k \geq 1} \mu(\Omega_k) \log \frac{1}{\mu(\Omega_k)} 2 \rho^{k+1} \mu(g^2) \\
& \leq &
2B^+ (\Phi) \sum_{k \geq 1} 2 \rho^{k+1} \mu(g^2) \alpha_{a_k}^+(\rho)
\end{eqnarray*}
where we used $\log (2 \rho^{k+1}) \leq 2\log (2 \rho^{k})$
for  $k \geq 1$ and the definition of $B^+ (\Phi)$.
Now consider the function $g_k=\ind_{[m,a_{k-1}[} +
\ind_{[a_{k-1},a_k[} \frac{g}{\sqrt{2\rho^{k-1}\mu(g^2)}} + \sqrt \rho
\ind_{[a_k,\infty)}$.  Since $g_k$ is non-decreasing, $g_k(m)=1$ and
$g_k(a_k)=\sqrt \rho$, we have
$$
\alpha_{a_k}^+ (\rho)\leq \int_m^{a_k} \Phi \left( \frac{g_k'}{g_k}
\right) g_k^2 d\nu \leq \frac{1}{2 \rho^{k-1} \mu(g^2)}
\int_{a_{k-1}}^{a_k} \Phi \left( \frac{g'}{g} \right) g^2 d\nu.
$$
Thus,
\begin{eqnarray*}
\int_{\Omega_1} g^2 \log \frac{g^2}{\mu(g^2)} d\mu
& \leq &
2 \rho^2 B^+ (\Phi) \sum_{k \geq 1}
\int_{a_{k-1}}^{a_k} \Phi \left( \frac{g'}{g} \right) g^2 d\nu \\
& \leq &
2 \rho^2 B^+ (\Phi) \int_{\Omega_0} \Phi \left( \frac{f'}{g} \right) g^2 d\nu\\
&\le & 2 \rho^2 B^+ (\Phi) \int_{[m,+\infty)} \Phi \left( \frac{f'}{f} \right) f^2 d\nu,
\end{eqnarray*}
where we have used that $f\le g$ and the monotonicity of $\Phi(t)/t^2$.

The third term in \eqref{eq:start} is estimated in a similar way. Finally one gets
\begin{eqnarray*}
   \ent_\mu(f^2) &\le& 16(1+\sqrt{2\rho})^2 C_P \int \left(\frac{f'}{f}\right)^2 f^2 d\mu+\\
  &&  2\rho^2 \max(B_+(\Phi),B_-(\Phi)) \int \Phi \left(\frac{f'}{f}\right) f^2 d\mu.
\end{eqnarray*}
Our hypotheses ensure that $H_\Phi(x) \ge \max (x^2, \Phi(x)/\Phi(1))$, hence
$$  \ent_\mu(f^2) \le \Big( 16(1+\sqrt{2\rho})^2 C_P +
    2\rho^2 \Phi(1)\max(B_+(\Phi),B_-(\Phi))\Big) \int H_\Phi \left(\frac{f'}{f}\right) f^2 d\mu.$$
\end{proof}

\begin{proposition} \label{prop:min}
   Let $\mu$ be a probability measure with median $m$ and $\nu$ a
   non-negative measure, on $\dR$. Assume that $\nu$ is absolutely continuous with respect to
   Lebesgue measure. Let $\Phi$ be a nice Young function
   and $H_\Phi$ its modification (see \eqref{eq:m}).

   Define the quantities $\alpha_x^+$ for $x>m$ and $\alpha_x^-$ for
   $x<m$ as follows
   $$
   \widetilde \alpha_x^+:= \inf \left\{ \int_m^x H_\Phi
      \left(\frac{f'}{f} \right) f^2 d\nu, f \mbox{non-decreasing},
      f(m)=0, f(x)=1 \right\},
   $$
   $$
   \widetilde \alpha_x^-:= \inf \left\{ \int_x^m H_\Phi
      \left(\frac{f'}{f} \right) f^2 d\nu , f \mbox{non-increasing},
      f(x)=1, f(m)=0 \right\} .
   $$
   Let
\begin{eqnarray*}
\widetilde B^+ & := &
\sup_{x >m} \mu([x,\infty)) \log \left( 1+  \frac{1}{2 \mu([x,\infty))} \right)
\frac{1}{\widetilde \alpha_x^+} , \\
\widetilde B^- & := &
\sup_{x < m} \mu((-\infty,x]) \log \left(1+ \frac{1}{2\mu((-\infty,x])} \right)
 \frac{1}{\widetilde \alpha_x^-} .
\end{eqnarray*}
If $C$ is a constant such that for any locally Lipschitz $f: \dR \rightarrow \dR$,
\begin{equation} \label{eq:mls2}
\ent_\mu (f^2) \leq C \int H_\Phi \left(\frac{f'}{f} \right) f^2 d\nu,
\end{equation}
then
$$C \geq \max(\widetilde B^+, \widetilde B^-).
$$

\end{proposition}

\begin{proof}
   Fix $x_0>m$ and consider a non-decreasing function $f$ with
   $f(m)=0$ and $f(x_0)=1$. Consider the function $\widetilde f = f
   \ind_{[m,x_0[} + \ind_{[x_0,\infty)}$.  Following
   \cite{bartr03sipm} and starting with the variational expression of
   entropy (see e.g. \cite[chapter 1]{Ane}),
\begin{eqnarray*}
\ent_\mu(\widetilde f^2)
& = &
\sup\left\{ \int \widetilde f^2 g \,d\mu , \int e^g d\mu \leq 1 \right\} \\
& \geq &
\sup\left\{ \int_{[m,+\infty)} \widetilde f^2 g \,d\mu , g \geq 0 \mbox{ and }
\int_{[m,+\infty)} e^g d\mu \leq 1 \right\} \\
&\ge &
\sup\left\{ \int_{[x_0,+\infty)}  g\, d\mu , g \geq 0 \mbox{ and }
\int_{[m,+\infty)} e^g d\mu \leq 1 \right\} \\
& = &
 \mu([x_0,\infty)) \log \left( 1+  \frac{1}{2 \mu([x_0,\infty))} \right)
\end{eqnarray*}
where the first inequality relies on the fact that $\widetilde f = 0$
on $(-\infty, 0]$ (hence the best is to take $g=-\infty$ on $(-\infty,
0]$). The latter equality follows from \cite[Lemma 6]{bartr03sipm}
which we recall below. Applying the modified logarithmic Sobolev
inequality to $\widetilde f$, we get
$$
\mu([x_0,\infty)) \log \left( 1+ \frac{1}{2 \mu([x_0,\infty))}
\right) \leq C \int_m^{x_0} H_\Phi \left(\frac{f'}{f} \right) f^2 d\nu
.
$$
Optimizing over all non-decreasing functions $f$ with $f(m)=0$ and
$f(x_0)=1$, we get
$$
\mu([x_0,\infty)) \log \left( 1+ \frac{1}{2 \mu([x_0,\infty))}
\right) \leq C \widetilde \alpha_{x_0}^+ .
$$
Hence $C\ge \widetilde B^+$.  A similar argument on the left of the
median yields $C\ge \widetilde B^-$.
\end{proof}

\begin{lemma}[\cite{bartr03sipm}]
   Let $Q$ be a finite measure on a space $X$. Let $K>Q(X)$ and let
   $A\subset X$ be measurable with $Q(A)>0$. Then
   $$
   \sup\left\{\int_X \ind_A h\, dQ;\; \int_X e^{h}dQ\le K
      \,\mathrm{and} \, h\ge 0 \right\} = Q(A)\log\left( 1+
      \frac{K-Q(X)}{Q(A)}\right).
   $$
\end{lemma}

\begin{remark}\label{rem:log}
   For $x \in (0,\frac 12)$, $\frac 34 \log \frac 1 x \leq \log( 1 +
   \frac{1}{2x} ) \leq \log \frac 1 x$.  Hence $B^+$ is comparable to
   $$
   \sup_{x >m} \mu([x,\infty)) \log
   \left(\frac{1}{\mu([x,\infty))}\right) \frac{1}{\widetilde
     \alpha_x^+}
   $$
   and similarly for $\widetilde B^-$.
\end{remark}

\bigskip

In order to turn the previous abstract results into efficient
criteria, we need more explicit estimates of the quantities
 $\alpha_x$ and $\widetilde \alpha_x$.


\section{The  example of power functions: $\Phi(x)=|x|^q$, $q \geq 2$.}

In this section we set $\Phi(x)=\Phi_q(x)=|x|^q$, with $q \geq 2$. Its modification is
 $H(x)=H_q(x)=\max(x^2,|x|^q)$. The constants $\alpha_x^\pm$
and $\widetilde \alpha_x^\pm$ are defined accordingly
as in Proposition~\ref{prop:maj} and Proposition~\ref{prop:min}.

The definition of $\alpha^+_x$ is simpler than the one of $\widetilde \alpha^+_x$.
Indeed it involves only $\Phi_q$. This allows the following easy estimate.

\begin{lemma} \label{lem:maj}
Assume that $\nu$ is absolutely continuous
with respect to the Lebesgue measure on $\dR$,
with density $n$.
Then for  $x>m$
$$
\frac{1}{\alpha_x^+}
\leq 2^{q-2}
\left(\int_m^x  n^\frac{-1}{q-1} \right)^{q-1} .
$$
\end{lemma}

\begin{proof}
Fix $x>m$. Let $q^*$ be such that
$\frac 1q + \frac{1}{q^*} =1$.
Consider a non-decreasing function $f$ with $f(m)=1$ and $f(x)=2$.
We assume without loss of generality that $\int_m^x |f'|^qf^{2-q} d\nu$ and $\int_m^xn^{-q^*/q}$ are finite.
By  Hölder's inequality (valid also when $n$ vanishes), we have
\begin{equation*}
 1 =  \int_m^x f' \le
\left( \int_m^x |f'|^q n \right)^\frac 1q
\left( \int_m^x n^{-\frac{q^*}{q}} \right)^\frac{1}{q^*}
\leq
\left(2^{q-2} \int\left| \frac{f'}{f}\right|^q f^2 d\nu \right)^\frac 1q
\left( \int_m^x n^{-\frac{q^*}{q}} \right)^\frac{1}{q^*},
\end{equation*}
where we used the bounds $f\le 2$ and $q\ge 2$. The result follows at once.
\end{proof}
A similar bound is available for $\alpha_x^-$ when $x<m$. Next we study the quantities
$\widetilde\alpha_x^+$. They are estimated by testing the inequality on specific functions,
as in the proofs of Hardy's inequality. However the presence of the modification $H_q$
creates complications, and we are lead to make additional assumptions.
We also omit the corresponding bound on $\widetilde \alpha_x^-$.

\begin{lemma} \label{lem:min}
Let $\nu$ be a non-negative measure absolutely continuous with
respect to the Lebesgue measure on $\dR$,
with density $n$.
Assume that there exists $\varepsilon>0$ such that for every $x>m$, it holds
\begin{equation}
   \label{eq:epsilon}
   (q-1) n(x)^{\frac{-1}{q-1}} \ge \varepsilon \int_{m}^x n(u)^{\frac{-1}{q-1}}du.
\end{equation}
Then for $x > m$, the quantity
$$
\widetilde \alpha_x^+=
\inf \left\{ \int_m^x H_q \left( \frac{f'}{f} \right) f^2 d\nu,
f \mbox{ non-decreasing}, f(m)=0, f(x)=1 \right\} .
$$
verifies
$$
\frac{1}{\widetilde \alpha_x^+} \geq
\frac{\min\big(\varepsilon^{q-2},1\big)}{ (q - 1)^{q -1}}
\left(\int_m^x  n(u)^{\frac{-1}{q-1}}du \right)^{q-1} .
$$
\end{lemma}

\begin{proof}
Fix  $x>m$. Then define
$$
f_x(t)= \left(
\frac{\displaystyle \int_m^t n^{\frac{-1}{q-1}}}
{\displaystyle \int_m^x  n^{\frac{-1}{q-1}}} \right)^{q-1}
\ind_{[m,x]} +  \ind_{(x,\infty)} .
$$
Note that $f_x$ is non-decreasing and satisfies
$f_x(m)=0$ and $f_x(x)=1$.
Thus,
$$
\widetilde \alpha_x^+
\leq
\int_m^x H_q \left( \frac{f_x'}{f_x} \right) f_x^2 d\nu .
$$
Furthermore \eqref{eq:epsilon} yields for $t\in (m,x)$,
$$
\frac{f_x'(t)}{f_x(t)} =
\frac{(q -1) n(t)^{\frac{-1}{q-1}}}{\displaystyle
\int_m^t  n^{\frac{-1}{q-1}} } \ge \varepsilon.
$$
Since
$H_q (t) \leq
\max \left( \frac{1}{\varepsilon^{q -2}},1 \right)
t^q$ for $t\in[ \varepsilon,\infty)$,
it follows, after some  computations, that
\begin{eqnarray*}
\int_m^x H_{q} \left( \frac{f_x'}{f_x} \right) f_x^2 d\nu
& \leq &
\max \left( \frac{1}{\varepsilon^{q-2}},1 \right)
\int_m^x  \left( \frac{f_x'}{f_x} \right)^q f_x^2 d\nu \\
& = &
\max \left( \frac{1}{\varepsilon^{q -2}},1 \right)
\frac{ (q - 1)^{q -1} }{\displaystyle
 \left(\int_m^x
n^{\frac{-1}{q-1}} \right)^{q -1}} \cdot
\end{eqnarray*}
This is the expected result.
\end{proof}

The next result provides a simple condition ensuring Hypothesis~\eqref{eq:epsilon} to hold
\begin{lemma}\label{lem:epsilon2}
   For a function $n(x)=e^{-V(x)}$ defined for $x\ge m$. Assume that for $x\in[m,m+K]$ one has $|V(x)|\le C$
   and that $V$ restricted to $[m+K,+\infty)$ is $C^1$ and verifies $V'(x)\ge \delta>0$, $x\ge m+K$.
   Then for $x\ge m$, one has
   $$ (q-1) n(x)^{\frac{-1}{q-1}} \ge \varepsilon \int_{m}^x n^{\frac{-1}{q-1}},$$
  where $\displaystyle\varepsilon=\frac{1}{\frac1\delta + \frac{K}{q-1} e^{2C/(q-1)}}>0$.
\end{lemma}
\begin{proof}
   Note that $V(x)\ge -C$ is actually valid for all $x\ge m$. If $x\le m+K$, simply write
   $$ \int_{m}^x n^{-\frac{1}{q-1}}=\int_{m}^x e^{\frac{V}{q-1}}\le Ke^{\frac{C}{q-1}} \le K e^{\frac{2C}{q-1}}
e^{\frac{V(x)}{q-1}}= K e^{\frac{2C}{q-1}} n(x)^{\frac{-1}{q-1}}.$$
If $x>m+K$, then
\begin{eqnarray*}
    \int_{m}^x e^{\frac{V}{q-1}} &\le&  Ke^{\frac{C}{q-1}}+ \int_{m+K}^x e^{\frac{V}{q-1}} \\
 &\le &  Ke^{\frac{2C}{q-1}} e^{\frac{V(x)}{q-1}}+ \frac{1}{\delta}\int_{m+K}^x V' e^{\frac{V}{q-1}} \\
 &=&  Ke^{\frac{2C}{q-1}} e^{\frac{V(x)}{q-1}}+ \frac{q-1}{\delta} \Big(e^{\frac{V(x)}{q-1}}- e^{\frac{V(m+K)}{q-1}}\Big)\\
  &\le & \left(Ke^{\frac{2C}{q-1}}+ \frac{q-1}{\delta} \right)e^{\frac{V(x)}{q-1}}.
\end{eqnarray*}
\end{proof}

\begin{theorem} \label{th:phialpha}
Let $\mu$ be a probability measure on $\dR$ with median $m$.
Let $\nu$ be a positive measure absolutely continuous with respect to the Lebesgue
measure with density $n$. Let $C_P\in(0,+\infty]$ be the optimal constant so that the Poincaré
 inequality \eqref{eq:poincare} holds.
Fix $q \geq 2$ and define
\begin{eqnarray*}
B^+_q &: = &
\sup_{x >m} \mu([x,\infty)) \log \frac{1}{\mu([x,\infty))}
\left(\int_m^x   n^{\frac{-1}{q-1}} \right)^{q -1}   , \\
B^-_q &: = &
\sup_{x <m} \mu((-\infty,x]) \log \frac{1}{\mu((-\infty,x])}
\left(\int_x^m  n^{\frac{-1}{q-1}} \right)^{q-1}.
\end{eqnarray*}
Let $\kappa_q\in (0,+\infty]$ be the best constant such that every locally Lipschitz
$f: \dR \rightarrow \dR$ satisfies
\begin{equation} \label{eq:phialpha}
\ent_\mu (f^2) \leq
\kappa_q \int H_q \left( \frac{f'}{f} \right) f^2 d\nu.
\end{equation}
Then $$ \displaystyle \kappa_q \leq  235 C_P + 2^{q+1} \max(B^+_q,B^-_q).$$
If there exists $\varepsilon>0$ such that for all $x\neq m$,
$$ (q-1) n(x)^{\frac{-1}{q-1}}\ge \varepsilon \int_{\min(x,m)}^{\max(x,m)} n^{\frac{-1}{q-1}},$$
then it is also true that
$$
\kappa_q \ge  \max\left(2C_P, \frac{3\min\big(\varepsilon^{q-2},1\big)}{4 (q-1)^{q-1}} \max\big(B_q^+,B_q^-\big)\right).
$$
\end{theorem}

\begin{proof}
The upper bound is immediate from Proposition~\ref{prop:maj} and Lemma~\ref{lem:maj} (and its obvious counterpart
on the left of the median).
The lower bound $\kappa_q\ge 2C_P$ is well known, see \cite{gentilgm}. It follows from applying the modified
log-Sobolev inequality to $f=1+tg$ where $g$ is a bounded function and $t$ goes to zero. Indeed $\ent_\mu((1+tg)^2)$
tends to $2\var_\mu(g)$ in this case. The lower bound in terms of $B_q^{\pm}$ is a direct consequence of
Proposition~\ref{prop:min}, Remark~\ref{rem:log} and Lemma~\ref{lem:min}.
\end{proof}

The following classical lemma (see e.g. \cite[Chapter 6]{Ane}) allows to estimate the integrals appearing in $B_q^\pm$.
\begin{lemma} \label{lem:equiv}
Let  $\Psi:[a,+\infty) \to \mathbb R^+$ be a locally bounded function. Assume that it is  ${\cal C}^2$ in a neighborhood of $+\infty$
and satisfies
$\liminf_\infty \Psi' >0 $.
\begin{enumerate}
\item If $\lim_\infty \Psi ''(x)/ \Psi '(x)^2 = 0$ then for $x$ growing to   infinity
$$
\int_a^x e^{\Psi (t)}dt \sim \frac{e^{\Psi (x)}}{\Psi '(x)} ,
\qquad \mbox{and}
\qquad
\int_x^{+\infty} e^{-\Psi (t)}dt \sim \frac{e^{-\Psi (x)}}{\Psi '(x)}.
$$

\item If for $x\ge x_0$  and $\varepsilon,A>0$, it holds $-1+\varepsilon \le \frac{\Psi''(x)}{\Psi'(x)^2}\le A$, then for $x\ge x_0$
$$ \frac{1}{1+A} \frac{e^{-\Psi (x)}}{\Psi '(x)} \le
\int_x^{+\infty} e^{-\Psi (t)}dt \le \frac{1}{\varepsilon}  \frac{e^{-\Psi (x)}}{\Psi '(x)}.$$
\end{enumerate}
\end{lemma}
As an application we obtain a workable criterion for satisfying a modified log-Sobolev inequality
with function $H_q$.

\begin{theorem}\label{th:critq}
   Let $q\ge 2$. Let $d\mu(x)=e^{-V(x)}dx$ be a probability measure on $ \mathbb R$.
Assume that $V:\mathbb R\to \mathbb R$ is locally bounded, and $C^2$ in neighborhoods of $+\infty$ and
   $-\infty$ with
\smallskip

(i) $\displaystyle\liminf_{|x|\to\infty} \mathrm{sign}(x) V'(x) >0$

(ii)  $\displaystyle \lim_{|x|\to \infty} \frac{V''(x)}{V'(x)^2}=0.$

\noindent
Then, there exists  $\kappa<+\infty$ such that for every locally Lipschitz $f$,
 $$ \ent_\mu(f^2) \le \kappa \int H_q\left(\frac{f'}{f}\right)f^2 d\mu$$
if and only if
$$ \limsup_{|x|\to \infty} \frac{V(x)}{|V'(x)|^q}<\infty.$$
\end{theorem}
\begin{remark}
   The condition on $V''/(V')^2$ can be relaxed to $-1<\liminf \frac{V''}{(V')^2} \le
  \limsup \frac{V''}{(V')^2}<\frac{1}{q}$. See Section~\ref{sec:general} where this is done in the general case.
\end{remark}
\begin{proof}
   Combining Theorem~\ref{th:hpoinc} (for $\nu=\mu$)  with Lemma~\ref{lem:equiv} shows
   that $\mu$ satisfies a Poincaré inequality. The hypotheses of Lemma~\ref{lem:epsilon2} are satisfied,
   therefore we may apply the two results in Theorem~\ref{th:phialpha}. It follows that
   $\mu$ satisfies the modified log-Sobolev inequality if and only if the quantities $B_q^+$ and $B_q^-$
   are finite. The potential $V$ being locally bounded we only have to care about large values of the
   variables. Applying Lemma~\ref{lem:equiv} again, we see that for $x$ large
    $$ \mu([x,+\infty)) \log\left(\frac{1}{\mu([x,+\infty))}\right) \left(\int_m^x e^{\frac{V}{q-1}}\right)^{q-1}
   \sim  \frac{V(x)+\log V'(x)}{V'(x)^q}\cdot$$
  Hence  $B_q^+$ is finite if and only if $\frac{V+\log V'}{(V')^q}$ has a finite upper limit at $+\infty$.
  By $(i)$, the term $V'$ is bounded away from 0 in the large. Thus $\log(V')/(V')^q$ is bounded and only
  $V/(V')^q$ matters.  A similar argument allows to deal with $B_q^-$.
\end{proof}

As a direct consequence we recover Theorem~3.1 of  Gentil, Guillin and Miclo \cite{gentilgm}.
\begin{corollary} \label{cor:mualpha}
Fix $q \geq 2$ and define its dual exponent $q^*$ by
$\frac{1}{q}+ \frac{1}{q^*} =1$. Let $p \geq 1$ and
 $d\mu_p(x)=Z_p^{-1}e^{-|x|^{p}}dx$.
Then there exists a constant $C_{p,q}<+\infty$ such that every locally Lipschitz
$f: \dR \rightarrow \dR$ satisfies
$$
\ent_{\mu_p} (f^2) \leq C_{p,q}
\int H_q\left( \frac{f'}{f} \right) f^2 d\mu_p
$$
if and only if $p\ge q^*$.
\end{corollary}

\begin{remark}
Bobkov and Ledoux \cite{bobkl97pitc} proved that a measure satisfies a Poincaré inequality if and only if it satisfies
a modified logarithmic Sobolev inequality with  function $H(t)=t^2\ind_{|t|\le t_0}$. This equivalence
 yields an improvement of the concentration inequalities that one
can deduce from a Poincaré inequality. It is natural to conjecture equivalences between general
modified log-Sobolev inequalities and inequalities involving $\int (f')^2d\mu$.
Under the hypotheses of the above theorem, Proposition~15 in \cite{bartr03sipm} shows that
the condition $ \limsup_{|x|\to \infty} \frac{V(x)}{|V'(x)|^q}<\infty$ is also equivalent
to $\mu$ satifying the following  Lata{\l}a-Oleszkiewicz inequality: there exists $\lambda<+\infty$ such that
for all locally Lipschitz $f$,
 $$ \sup_{\theta\in [1,2)} \frac{\int f^2 d\mu-\left(\int |f|^\theta d\mu\right)^{2/\theta}}{(2-\theta)^{2/q}} \le
   \lambda \int (f')^2 d\mu.$$
Hence, under the hypotheses of Theorem~\ref{th:critq}, a measure satisfies the latter inequality if and only if
it satisfies a modified log-Sobolev inequality with function $H_q$.

\end{remark}

\begin{remark}
   It is known that general modified log-Sobolev inequalities imply
so-called transportation cost inequalities, see \cite{bobkgl01hhje}. Criteria
 for measures on the line to satisfy such inequalities  have been obtained
recently by Gozlan \cite{gozl06ctti}, after a breakthrough of Cattiaux and Guillin \cite{cattg06qtci}. It is interesting
to compare his result with  Theorem~\ref{th:critq}.
\end{remark}


\section{More general cases} \label{sec:general}
The results of the previous section extend to more general functions $\Phi$.
Now, we show  how to reach them.
In order to obtain workable versions of Propositions~\ref{prop:maj} and \ref{prop:min},
we need explicit lower bounds on
 $\alpha^+_x$ and  $\alpha^-_x$ as well as upper bounds on $\widetilde \alpha_x^+$ and $\widetilde\alpha_x^-$.
Actually our methods also allow bounds in the other direction, but we omit them as they have
no other use than showing that the bounds are rather good. By symmetry we shall discuss only
$\alpha_x^+$ and $\widetilde \alpha_x^+ $.

In all this section, $\Phi$ stands for a nice Young function,
$\Phi^*$ for its conjugate and
$\nu$ for a non-negative measure on $\dR$.

\subsection{Lower bounds on  $\alpha_x$. Sufficient conditions}
Given   $x>m$, we have set
$$
\alpha_x^+= \inf \left\{ \int_m^x \Phi \left( \frac{f'}{f} \right) f^2 d\nu,
f \mbox{ non-decreasing}, f(m)=1, f(x)=2 \right\} .
$$
The following simple lower bound is available
\begin{eqnarray*}
   \alpha_x^+&\ge &  \inf \left\{ \int_m^x \Phi \left( \frac{f'}{2} \right)  d\nu,
f \mbox{ non-decreasing}, f(m)=1, f(x)=2 \right\} \\
   &\ge &  \inf \left\{ \int_m^x \Phi \left( \frac{g}{2} \right)  d\nu,
  g\ge 0, \int_m^x g(u) \, du =1 \right\} = \beta_x\Big(\frac{1}{2}\Big),
\end{eqnarray*}
where we have set for $a>0$,
$$
\beta_x(a):=\inf\left\{ \int_m^x \Phi (g)\, d\nu\,;\; g\ge 0 \;\mathrm{and} \,
\int_m^x g(t)\,dt=a \right\} .
$$

The infimum is  evaluated in the next lemma.
A similar result has been recently established by Arnaud Gloter \cite{gloter}.
The statement involves the following new notation.
 The left inverse of a non-decreasing function $f$ is defined by
$f^{-1}(x) :=\inf\{y;\; f(y)\ge u\}.$ Also for a non-decreasing function $\Psi$ on $\mathbb R^+$ with limits $0$
at $0$ and $+\infty$ at $+\infty$ but not necessarily convex, we define for a measurable function on $\mathbb R$
$$ \left\|g \right\|_{\Psi}
:=\inf \left\{\delta >0;\; \int_{\mathbb R} \Psi\left(\frac{|g|}{\delta}\right)\le 1\right\},$$
which needs not be a norm.
\begin{lemma} \label{prop:op}
Assume that $\nu$ is absolutely continuous
with respect to the Lebesgue measure on $\dR$,
with density $n$.
Then,
$$
\beta_x(a) \ge \int_m^x \Phi \left( {\Phi_r'}^{-1} \left( \frac{\gamma_{x,a}}{n } \right) \right) d\nu
$$
where
 $$\gamma_{x,a}:= \sup\left\{\lambda\ge 0;\;
\int_m^x {\Phi'_r}^{-1} \left( \frac{\lambda}{n(u)} \right) du \le a
 \right\}=\left(\left\| \frac{\ind_{[m,x]}}{n}\right\|_{\frac1a {\Phi'_r}^{-1}}\right)^{-1},$$
and ${\Phi'_r}^{-1}$ is the left inverse of the right derivative of $\Phi$.

Moreover, if $\Phi'_r$ is strictly increasing and satisfies the following doubling  condition:  there
exists $K>1$ such that for all $x\ge 0$,  $\Phi'_r(Kx)\ge 2\Phi'_r(x)$, then
when $\gamma_{x,a}\neq 0$,
  $$ \int_m^x {\Phi'_r}^{-1} \left( \frac{\gamma_{x,a}}{n(u)} \right) du = a \quad\mbox{and}\quad
     \beta_x(a) = \int_m^x \Phi \left( {\Phi_r'}^{-1} \left( \frac{\gamma_{x,a}}{n } \right) \right) d\nu.$$
 \end{lemma}

\begin{proof}
If the set of points in $[m,x]$ where $n$ vanishes has positive Lebesgue measure, it is plain
that $\beta_x(a)=\gamma_{x,a}=0$ and the claimed result is obvious. Hence we may assume that
almost every  $t\in [m,x]$ verifies $n(t)>0$. We also assume that $\gamma_{x,a}>0$ otherwise there
is nothing to prove. Let us start with a nonnegative  function $g$ on $[m,x]$ with $\int_m^x g=a$ and
$\int_m^x \Phi(g)\, d\nu<\infty$. For $\lambda>0$, and almost every $t\in [m,x]$, $n(t)\neq 0$ and
 Young's inequality yields
$$ g(t) \le \frac{n(t)}{\lambda} \left( \Phi(g(t))+\Phi^*\Big(\frac{\lambda}{n(t)}\Big)\right),$$
where $\Phi^*(u) :=\sup_{y \geq 0} \{uy - \phi(y)\}$.
The analysis of equality cases in Young's inequality leads us to introduce
$$g_\lambda(t):=\inf\left\{x\ge 0;\; \Phi'_r(x)\ge \frac{\lambda}{n(t)} \right\}={\Phi'_r}^{-1}\Big( \frac{\lambda}{n(t)}
  \Big).$$
Since $\Phi'_r$ is right continuous and vanishes at 0, one has  $\Phi'_r(g_\lambda(t))
   \ge \frac{\lambda}{n(t)}\ge \Phi'_\ell(g_\lambda(t))$  (at least when $n(t)\neq 0$). By convexity this yields
 $$ \Phi^*\Big(\frac{\lambda}{n(t)}\Big)=\sup_{y\ge 0} \left\{\frac{\lambda}{n(t)}y-\Phi(y)\right\}=
       \frac{\lambda}{n(t)}g_\lambda(t)-\Phi\big(g_\lambda(t)\big).$$
Combining this  with the latter inequality gives
 $$ n(t) \Phi(g(t)) \ge n(t) \Phi(g_\lambda(t)) + \lambda(g(t)-g_\lambda(t)).$$
If $\lambda$ is chosen so that $\int_m^x g_\lambda \le a$, integrating the previous relation on $[m,x]$
implies that $\int_m^x \Phi(g) \, d\nu \ge \int_m^x \Phi(g_\lambda) \, d\nu$. Optimizing on $g$ and $\lambda$ satisfying the
above conditions, we obtain
$$ \beta_x(a)\ge \sup\left\{ \int_{m}^x \Phi\left({\Phi'_r}^{-1}\Big(\frac{\lambda}{n}\Big)\right) \, d\nu\right\},$$
where the supremum is taken above all $\lambda$ with $\int_m^x {\Phi'_r}^{-1}(\frac{\lambda}{n}) \le a$.
By definition $\gamma_{x,a}$ is the supremum of such $\lambda$'s. Using that a left inverse is left continuous,
we conclude that
$$ \beta_x(a)\ge  \int_{m}^x \Phi\left({\Phi'_r}^{-1}\Big(\frac{\gamma_{x,a}}{n}\Big)\right) \, d\nu.$$

If we also know that $\Phi'_r$ is strictly increasing, then its left inverse is continuous.
Moreover the doubling condition: $2\Phi'_r(x)\le \Phi'_r(Kx)$ translates to the left inverse as a so-called $\Delta_2$
 condition: for all $x\ge 0$, ${\Phi'_r}^{-1}(2x)\le K{\Phi'_r}^{-1}(x)$.
Hence for every positive real numbers $\lambda_1<\lambda_2$ and every $x\ge 0$,
  $$ {\Phi'_r}^{-1}(\lambda_1x)\le{\Phi'_r}^{-1}(\lambda_2 x)\le
      {\Phi'_r}^{-1}\left(2^{\left\lceil\frac{\log(\lambda_2/\lambda_1)}{\log2} \right\rceil} \lambda_1 x\right)
    \le K^{\left\lceil\frac{\log(\lambda_2/\lambda_1)}{\log2} \right\rceil}  {\Phi'_r}^{-1}(\lambda_1x).$$
Consequently the family of integrals $\left(\int_{m}^x  {\Phi'_r}^{-1}(\frac{\lambda}{n})\right)_{\lambda>0}$
are either simultaneously infinite or simultaneously finite. In the former situation one gets $\gamma_{x,a}=0$
whereas in the latter, the function $\lambda\mapsto \int_{m}^x  {\Phi'_r}^{-1}(\frac{\lambda}{n})$ is continuous
by dominated convergence and varies from 0 to $+\infty$ (recall that we reduced to $n>0$ almost everywhere on $[m,x]$).
Hence it achives the value $a>0$ for at least one $\lambda$ and the smallest of them is $\gamma_{x,a}$.
The function $g:=\frac{\gamma_{x,a}}{n}$ satisfies $\int_m^x g=a$ and
$$\int_m^x \Phi(g)\, d\nu=\int_m^x \Phi \left( {\Phi_r'}^{-1} \left( \frac{\gamma_{x,a}}{n } \right) \right) d\nu.$$
Hence the latter quantity coincides with $\beta_x(a)$.
\end{proof}

Under natural assumptions on the rate of growth of $\Phi$ we obtain a simpler bound on $\beta_x(a)$.
\begin{proposition} \label{prop:beta2}
Assume that $\nu$ is absolutely continuous
with respect to the Lebesgue measure on $\dR$,
with density $n$.  Assume that $\Phi$ is a strictly convex nice Young function such that
on $\mathbb R^+$ the function $\Phi(x)/x^2$ is non-decreasing and the function $\Phi(x)/x^\theta$ is non-increasing, where $\theta >2$. Then for all $a>0$,
  $$ \beta_x(a)\ge \frac{a\, \gamma_{x,a}}{\theta}\cdot$$
\end{proposition}

\begin{proof}
Assume as we may that $\gamma_{x,a}>0$. We check that the hypothesis of the stronger part of the previous
lemma are satisfied. The strict convexity of $\Phi$ ensures that $\Phi'_r$ is strictly increasing.
It remains to check the doubling condition for this function. By differentiation, the monotonicity of
 $\Phi(x)/x^2$ and $\Phi(x)/x^\theta$ yields for $x\ge 0$,
   $$ 2\Phi(x) \le x \Phi'_r(x) \le \theta \Phi(x).$$
Combining  these inequalities with the monotonicity of $\Phi(x)/x^2$ yields
 $$ \Phi'_r(\theta y) \ge 2 \frac{\Phi(\theta y)}{\theta y} \ge 2\theta \frac{\Phi(y)}{y} \ge 2\Phi'_r(y),$$
as needed. Applying the previous lemma, we obtain that $a=\int_{m}^x {\Phi'_r}^{-1}(\frac{\gamma_{x,a}}{n})$,
and
\begin{eqnarray*}
\beta_x(a)
& = &
\int_m^x \Phi \left( {\Phi'_r}^{-1} \left( \frac{\gamma_{x,a}}{n } \right) \right) d\nu\\
& \ge &
 \frac{1}{\theta} \int_m^x {\Phi'_r}^{-1} \left( \frac{\gamma_{x,a}}{n } \right)
\Phi'_r\left({\Phi'_r}^{-1} \left( \frac{\gamma_{x,a}}{n } \right)\right) n \\
& \ge & \frac{\gamma_{x,a}}{\theta} \int_m^x {\Phi'}^{-1} \left( \frac{\gamma_{x,a}}{n } \right)
 =\frac{a \,\gamma_{x,a} }{\theta},
\end{eqnarray*}
where we have used $F\big(F^{-1}(u)\big)\ge u$, valid  for any right-continuous function $F$.
\end{proof}

\begin{remark}
   When $\Phi(x)=|x|^q$, $\gamma_{x,a}$ and $\beta_x(a)$ are multiples of $(\int_m^x n^{\frac{-1}{q-1}})^{q-1}$. This is
consistent with Lemma~\ref{lem:maj}.
\end{remark}

Combining the Proposition~\ref{prop:maj} with the observation that $\alpha_x^+\ge \beta_x(1/2)$ and Proposition~\ref{prop:beta2},
we obtain the following criterion:
\begin{theorem}\label{th:critphi}
Let $\theta\ge 2$. Let $\Phi$ be a strictly convex nice Young function such that $\frac{\Phi(x)}{x^2}$ is non-decreasing
 and $\frac{\Phi(x)}{x^\theta}$ is non-increasing. Let $\mu$ be a probability measure on $\mathbb R$ with median $m$,
and let $d\nu(x)=n(x)\, dx$ be a measure on $\mathbb R$. Assume that they satisfy a Poincaré inequality \eqref{eq:poincare}
 with constant
$C_P$. Then for every locally Lipschitz function $f$ on $\mathbb R$, the following modified log-Sobolev inequality holds:
$$ \ent_\mu(f^2)\le \Big(235 C_P+16\,\theta\,\Phi(1) \max\big(C_-(\Phi),C_+(\Phi)\big)\Big) \int_{\mathbb R} H_\Phi\left(
  \frac{f'}{f}\right) f^2 d\nu,$$
with
\begin{eqnarray*}
   C_+(\Phi)&:=& \sup_{x>m} \mu\big([x,+\infty)\big) \log\Big(\frac{1}{\mu\big([x,+\infty)\big)}\Big)
                  \left\|\frac{\ind_{[m,x]}}{n} \right\|_{2 {\Phi'_r}^{-1}},\\
  C_-(\Phi)&:=& \sup_{x<m} \mu\big((-\infty,x]\big) \log\Big(\frac{1}{\mu\big((-\infty,x]\big)}\Big)
                  \left\|\frac{\ind_{[x,m]}}{n} \right\|_{2 {\Phi'_r}^{-1}}.
\end{eqnarray*}
\end{theorem}

\begin{lemma}\label{lem:maj-norme}
Let $\Phi$ be a differentiable, strictly convex nice Young function.
Assume that there exists $\theta>1$ such that $\Phi(x)/x^\theta$ is non-increasing
on $\mathbb R^+$.
Let $V:[m,+\infty)\to \mathbb R$ such that for all $x\in[m,m+K]$, it holds $|V(x)|\le C$.
Also assume that
$V$ is $C^2$ on  $[m+K, +\infty )$ and verifies for $x\ge m+K$,
 $$V'(x)>0 \quad \mbox{and} \quad \frac{V''(x)}{V'(x)^2}\le \frac{1}{\theta}.$$
Then for $\in ]m,m+K]$, it holds
 $\displaystyle  \left\|\frac{\ind_{[m,x]}}{e^{-V}} \right\|_{2{\Phi'}^{-1}} \le \frac{e^C}{\Phi'\Big(\frac{1}{4K}\Big)}$ and for
all $x>m+K$,
   $$ \left\|\frac{\ind_{[m,x]}}{e^{-V}} \right\|_{2{\Phi'}^{-1}} \le \max\left( \frac{e^C}{\Phi'\Big(\frac{1}{4K}\Big)},
           \frac{e^{V(x)}}{\Phi'\Big(\frac{V'(x)}{4\theta(\theta-1)}\Big)}\right).$$
\end{lemma}

\begin{proof}
Our hypotheses ensure that $\Phi'$ is a bijection of $[0;+\infty)$; its inverse is ${\Phi^*}'$.
In order to show that $\| f\|_\Psi\le \lambda$ it is enough to prove that $\int \Psi(|f|/\lambda) \le 1$.
Hence our task is to find $\varepsilon>0$ with  $\int_m^x 2{\Phi'}^{-1}(\varepsilon e^V)\le 1$.
We deal with the case $x\ge m+K$ (the remaining case is simpler and actually contained in the beginning
of the following argument):
 \begin{eqnarray*}
  \int_m^x {\Phi'}^{-1}(\varepsilon e^V) &=& \int_m^{m+K}  {\Phi'}^{-1}\big(\varepsilon e^{V(t)}\big)dt+
                      \int_{m+K}^x {\Phi'}^{-1}\big(\varepsilon e^{V(t)}\big)dt \\
              &\le & K {\Phi'}^{-1} (\varepsilon e^C)+   \int_{m+K}^x {\Phi^*}^{'}\big(\varepsilon e^{V(t)}\big)dt.       \end{eqnarray*}
The first term in the above sum is less than $1/4$ as soon as $\varepsilon \le e^{-C} \Phi'\Big(\frac{1}{4K}\Big).$
The last term is estimated by integration by parts:
\begin{eqnarray*}
      && \int_{m+K}^x {\Phi^*}'\big(\varepsilon e^{V(t)}\big)dt \,=\,
           \int_{m+K}^x  \varepsilon V'(t) e^{V(t)} {\Phi^*}'\big(\varepsilon e^{V(t)}\big)
                      \frac{1}{ \varepsilon V'(t) e^{V(t)}}\, dt \\
              &=& \frac{\Phi^*\big(\varepsilon e^{V(x)} \big)}{\varepsilon e^{V(x)} V'(x)} -
                  \frac{\Phi^*\big(\varepsilon e^{V(m+K)} \big)}{\varepsilon e^{V(m+K)} V'(m+K)}+
                   \int_{m+K}^x \frac{ \Phi^*\big(\varepsilon e^{V(t)}\big) }{\varepsilon e^{V(t)}}
                   \left(1+\frac{V"(t)}{V'(t)^2} \right) \, dt \\
              &\le&      \frac{\Phi^*\big(\varepsilon e^{V(x)} \big)}{\varepsilon e^{V(x)} V'(x)}
                          +\Big(1+\frac{1}{\theta} \Big) \int_{m+K}^x
                          \frac{ \Phi^*\big(\varepsilon e^{V(t)}\big) }{\varepsilon e^{V(t)}}  dt \\
               &\le&      \frac{\Phi^*\big(\varepsilon e^{V(x)} \big)}{\varepsilon e^{V(x)} V'(x)}
                          +\Big(1-\frac{1}{\theta^2} \Big) \int_{m+K}^x
                           {\Phi^*}'\big(\varepsilon e^{V(t)}\big) \, dt,
\end{eqnarray*}
where we have used in the last line the inequality
$\Phi^*(x) \le \left(1- \frac{1}{\theta} \right) x {\Phi^*}'(x)$, which
follows from our hypotheses by Lemma~\ref{lem:a1}. The term $\int  {\Phi^*}'\big(\varepsilon e^V\big)$  appears
on both sides of the inequality. So after rearrangement we get
 $$ \int_{m+K}^x {\Phi^*}'\big(\varepsilon e^{V(t)}\big)dt
    \le \theta^2 \frac{\Phi^*\big(\varepsilon e^{V(x)} \big)}{\varepsilon e^{V(x)} V'(x)}
    \le \theta(\theta-1) \frac{{\Phi^*}'\big(\varepsilon e^{V(x)} \big)}{V'(x)}=
    \theta(\theta-1) \frac{\Phi'^{-1}\big(\varepsilon e^{V(x)} \big)}{V'(x)}\cdot$$
 Hence    $ \int_{m+K}^x {\Phi^*}'\big(\varepsilon e^{V(t)}\big)dt \le 1/4$ holds when
  $$ \varepsilon \le e^{-V(x)} \Phi'\left(\frac{V'(x)}{4\theta(\theta-1)}  \right).$$
  Finally for
  $$\varepsilon_0:= \min\left( e^{-C} \Phi'\Big(\frac{1}{4K}\Big),e^{-V(x)} \Phi'\left(\frac{V'(x)}{4\theta(\theta-1)}  \right)\right),$$
  we have shown that $ \int_m^x 2{\Phi'}^{-1}\big(\varepsilon_0 e^V\big) \le 1$. This concludes the proof.
 \end{proof}
Lemma~\ref{lem:maj-norme} allows to get more explicit versions of Theorem~\ref{th:critphi}. Here is an example
\begin{theorem}\label{th:critphi-simple}
 Let $\Phi$ be a strictly convex differentiable nice Young function on $\mathbb R^+$. Assume that
$\Phi(x)/x^2$ is non-decreasing and that there exists $\theta>2$ such that $\Phi(x)/x ^\theta$ is non-increasing.

Let $d\mu(x)=e^{-V(x)}dx$ be a probability measure on $\mathbb R$. Assume that $V$ is locally bounded, of class $\mathcal C^2$
in neighborhoods of $+\infty $ and $-\infty $ such that:
\begin{enumerate}
\item $\displaystyle\liminf_{|x|\to +\infty } \mathrm{sign}(x) V'(x)>0$,
\item $\displaystyle   -1< \liminf_{|x|\to +\infty } \frac{V''(x)}{V'(x)^2} \le \limsup_{|x|\to +\infty } \frac{V''(x)}{V'(x)^2}
   <\frac{1}{\theta}$,
\item $\displaystyle\limsup_{|x|\to +\infty }\frac{V(x)}{\Phi\big(|V'(x)| \big)}<+\infty$ .
\end{enumerate}
Then there exists a constant $\kappa<+\infty $ such that for all locally Lipschitz $f$ on $\mathbb R$
$$\ent_{\mu} (f^2) \le \kappa \int_{\mathbb R}H_\Phi\left(\frac{f'}{f}\right)f^2 d\mu.$$
\end{theorem}
\begin{proof}
   Combining hypothesis $(i)$ with Theorem~\ref{th:hpoinc} for $\nu=\mu$ and Lemma~\ref{lem:equiv} shows that
$\mu$ satisfies a Poincaré inequality. Our task is therefore to show that the numbers $C_+(\Phi),C_-(\Phi)$ in
the statement of Theorem~\ref{th:critphi} are finite. By symmetry we only deal with $C_+(\Phi)$. Since $V$ is
locally bounded and $t\log(1/t)$ is upper bounded on $(0,1]$, Lemma~\ref{lem:maj-norme} allows us to reduce the
problem to the finiteness of the upper limit when $x\to +\infty$ of
$$ \mu\big([x,+\infty )\big) \log\Big(\frac{1}{\mu\big([x,+\infty )\big)} \Big)
 \frac{e^{V(x)}}{\Phi'\left(\frac{ V'(x)}{4\theta(\theta-1)}\right)} \cdot$$
For shortness we set $T:=4\theta(\theta-1)>1$. Our assumptions imply that there exists $\varepsilon>$ such that
  for $x$ large enough $1\ge V''(x)/V'(x)^2\ge -1+\varepsilon$. Thus, the second part of Lemma~\ref{lem:equiv} shows that the
above  quantity is at most
 $$\frac{V(x)+\log\big(2V'(x) \big)}{\varepsilon V'(x)\Phi'\left(\frac{ V'(x)}{T}\right) } \le
  \frac{V(x)+\log\big(2V'(x) \big)}{\varepsilon T\Phi\left(\frac{ V'(x)}{T}\right) }
    \le T^{\theta-1} \frac{V(x)+\log\big(2V'(x) \big)}{\varepsilon \Phi\big(V'(x) \big)},$$
where we have used that $\Phi(x)/x^\theta$ is non-increasing.
Finally since $V'(x)$ is bounded below by a positive number for large $x$, the ratio of $\log V'$ to $\Phi(V')$ is upper
bounded in the large. Condition $(iii)$ allows to conclude.
\end{proof}

As a direct consequence we recover  the result by Gentil-Guillin and Miclo \cite{gentilgm2}
with slightly different conditions.
\begin{corollary} \label{cor:gen}
Let $\Psi$ be an even convex function on $\mathbb R$ such that $d\mu_{\Psi}(x)=e^{-\Psi(x)} dx$ is a probability
measure. Let $\alpha\in (1,2]$.
 Assume that for $x\ge x_0$, $\Psi$   is of class $\mathcal C^2$ with  $\Psi(x)/x^2$  non-increasing
  and   $\Psi(x)/x^\alpha$  non-decreasing, and that  $\limsup_\infty \frac{{\Psi}''}{{\Psi'}^2} < 1-\frac{1}{\alpha}$.

Then there exists
$C,D \in (0,+\infty)$ such that, setting $\mathcal H(x)= C \big( x^2 \ind_{|x|<D}+ \Psi^*(|x|) \ind_{|x|\ge D}\big)$,
 every locally Lipschitz
$f : \dR \rightarrow \dR$ verifies
$$
\ent_{\mu_\Psi}(f^2) \leq  \int_{\mathbb R}\mathcal H \left( \frac{f'}{f} \right) f^2 d\mu_\Psi .
$$
\end{corollary}

\begin{remark} \label{rem:gen}
If for some $\varepsilon\in (0,1)$, $\Psi^\varepsilon$ is concave in the large, then
$\lim_\infty \frac{{\Psi}''}{{{\Psi}'}^2} = 0$.
\end{remark}

\begin{proof}
   We apply Theorem~\ref{th:critphi-simple} with a suitable function $\Phi$.
 We choose $x_1>x_0$ such that $\Psi(x_1)>1$ and $\Psi'(x_1)>1$.
 Our monotonicity assumptions ensure that for $x\ge x_0$, $\alpha \Psi(x)\le x\Psi'(x)\le 2\Psi(x)$.
Let $\beta=\frac{x_1\Psi'(x_1)}{\Psi(x_1)}\in [\alpha,2]$, and set for $x\ge 0$
 $$ f(x)= \Psi(x_1) \left(\frac{x}{x_1}\right)^\beta \ind_{x<x_1}+ \Psi(x) \ind_{x\ge x_1}.$$
One easily checks that $f$ is convex of class $C^1$, and that on $\mathbb R^+$, $f(x)/x^\alpha$ is non-decreasing
whereas $f(x)/x^2$ is non-increasing. By Lemma~\ref{lem:a1} the conjugate function is such that $f^*(x)/x^2$ is
non-decreasing and $f^*(x)/x^{\alpha^*}$ is non-increasing for $x>0$ and $\alpha^*=\alpha/(\alpha-1)\ge 2$.
One easily checks that for a suitable constant $b$ and for $x\ge 0$
$$
f^*(x)=b x^{\beta^*}\ind_{x< \Psi'(1)} + \Psi^*(x) \ind_{x\ge \Psi'(1)}.
$$
Finally we set $\Phi(x)=f^*(x)+x^2$ in order to have a strictly
convex function with the same monotonicity properties, to which
Theorem~\ref{th:critphi-simple} may be applied for $V=\Psi$. Note
that obviously $\lim_{+\infty } \Psi'=+\infty$. Our assumptions
imply that $0\le \liminf \frac{\Psi"}{{\Psi'}^2}\le \limsup
\frac{\Psi"}{{\Psi'}^2}
 < 1-\frac{1}{\alpha}=\frac{1}{\alpha^*}$.
Our task is to show the boundedness of the upper limit at $+\infty $ of $\frac{\Psi}{\Phi(\Psi')}$.
For $x$ large enough,
$$ \frac{\Psi(x)}{\Phi(\Psi'(x))} \le \frac{\Psi(x)}{\Psi^*(\Psi'(x))} \le \frac{\alpha^* \Psi(x)}{ \Psi'(x) {\Psi^*}'(\Psi'(x))}
  = \frac{\alpha^* \Psi(x)}{ \Psi'(x) x} \le  \frac{\alpha^*}{\alpha},
$$
where we have used, in differential form, the fact that in the large $\Psi^*(x)/x^{\alpha^*}=f^*(x)/x^{\alpha^*}$ is non-increasing
and $\Psi(x)/x^\alpha$ is non-decreasing.
Since $\Psi $ is even,   Theorem~\ref{th:critphi-simple} ensures that
the measure $\mu_\Psi$ satisfies a modified log-Sobolev inequality with function
$H_\Phi$. One easily checks that for suitable choice of $C,D$, this function $H_\Phi$ is upper-bounded by the function
$\mathcal H$ of the claim.

\end{proof}

\subsection{Upper bounds on  $\widetilde \alpha_x$. Necessary conditions.}

Recall that we have set  for $x>m$,
$$
\widetilde \alpha_x^+= \inf \left\{
\int_m^x H_\Phi \left( \frac{f'}{f} \right) f^2 d\nu,
f \mbox{ non-decreasing}, f(m)=0, f(x)=1 \right\} ,
$$
where $H_\Phi$ stands for the modification of $\Phi$ (see \eqref{eq:m}).
In order to get necessary conditions for modified log-Sobolev inequalities to hold, we need upper
bounds on $\widetilde \alpha_x^+$. The next result provides an asymptotic estimate. Noting that $\widetilde \alpha_x^+ \geq \frac{1}{\Phi(1)} \alpha_x^+$  holds when $\Phi(x)/x^2$ is non-decreasing and comparing with the lower
bound on $\alpha_x^+$ given (in different notation) in Lemma~\ref{lem:maj-norme} shows that the bound is
of the right order.

\begin{proposition}\label{prop:al}
Let $\Phi$ be a twice differentiable, strictly convex, nice Young function. Assume that on $\mathbb R^+$
the function $\Phi(x)/x^2$ is non-decreasing, the functions $\Phi(x)/x^\theta$ and $\Phi'(x)/x^\eta$
are non increasing for some $\theta,\eta>0$. Also assume that there exists $\Gamma\in\mathbb R$ such that
for all $x,y\ge 0$, $\Gamma \Phi(xy)\ge \Phi(x)\Phi(y)$.

Let  $d\nu(x)=e^{-V(x)} \, dx$ be a measure on $\mathbb R$.
Assume furthermore that $V$ is $C^2$ in a neighborhood of $+\infty$, with
\begin{enumerate}
\item
$\displaystyle \liminf_{x\to +\infty} V'(x)>0$,
\item
$\displaystyle -1< \liminf_{x\to +\infty } \frac{V''(x)}{V'(x)^2} \le
\limsup_{x\to +\infty } \frac{V''(x)}{V'(x)^2}<\frac{1}{\max(\theta,\eta)}$.
\end{enumerate}
Then there exists a number $K$ depending only on $V$ and $\Phi$ such that for $x$ large enough,
$$
\widetilde \alpha_x^+
\leq
K e^{- V(x)} \Phi' \big(V'(x)\big).
$$
\end{proposition}

\begin{proof}
We shall prove the above inequality for $x\ge x_1>x_0>m$ where $x_0,x_1$ are  large enough.  We start with $\varepsilon\in (0,\liminf V')$  small enough to
have $\limsup \frac{V"}{{V'}^2} < \frac{1-\varepsilon}{\eta}$. We choose $x_0$
large enough to ensure that for $x \ge x_0$,
$$V'(x)>\varepsilon \mbox{ and }  -1\le \frac{V"(x)}{{V'(x)}^2}\le \min\left(\frac{1}{\theta}, \frac{1-\varepsilon}{\eta}\right).$$
 For $x\ge x_0$, let
$$
f_x(t):= \ind_{[x_0,x]} \,\int_{x_0}^t {\Phi'}^{-1} \left( c_x e^{V(u)} \right) du
 +  \ind_{]x,\infty)}
$$
where
$c_x >0$ is such that
$\int_{x_0}^x {\Phi'}^{-1} \left( c_x e^{V(u)} \right) du = 1$.
We also define $g_x := \Phi'(f_x)$.

The hypothesis on $\Phi'$ is equivalent to $t\Phi''(t)\le \eta\Phi'(t)$. Hence for $x>x_0$
$$
\frac{g_x'}{g_x}  = \frac{\Phi''(f_x) f_x'}{\Phi'(f_x)} \leq \eta\frac{f_x'}{f_x}.
$$
Since $g_x$ is non-decreasing and satisfies $g_x(m)=0$,
$g_x(x)=\Phi'(1)$, it follows that
$$
\widetilde{\alpha}_x^+ \le \int_m^x  H_\Phi \left( \frac{g_x'}{g_x} \right)
\left(\frac{g_x}{\Phi'(1)}\right)^2 d\nu
 \leq \Phi'(1)^{-2}
\int_{x_0}^x  H_\Phi
\left(\eta \frac{f_x'}{f_x} \right) \Phi'(f_x)^2 d\nu .
$$
 Lemma~\ref{lem:a3} ensures that the hypothesis  $\Gamma \Phi(xy)\ge \Phi(x)\Phi(y)$ transfers to
 ${\Phi'}^{-1} = {{\Phi^*}'}$. More precisely there exists another constant $\Gamma'$ such that
for $x,y\ge 0$, ${\Phi^*}'(xy) \le \Gamma' {\Phi^*}'(x){\Phi^*}'(y)$. Hence for $t\in(x_0,x)$,
\begin{eqnarray*}
\frac{f_x(t)}{f_x'(t)}
 =
\int_{x_0}^t
\frac{{\Phi^*}' \left( c_x e^{V(u)} \right)}
{{\Phi^*}' \left( c_x e^{V(t)} \right)} du
 \leq
 \Gamma' \int_{x_0}^t {\Phi^*}' \left( e^{V(u) - V(t)} \right) du .
\end{eqnarray*}
Using our assumption that $V'(x)>\varepsilon$ for $x\ge x_0$, and the inequality
${{\Phi^*}'}(x) \leq 2 x^{1/(\theta-1)} \Phi^*(1)$, for $x\in [0,1]$, a consequence of
Lemma~\ref{lem:a2} of the Appendix,
we obtain
$$
\frac{f_x(t)}{f_x'(t)}
\leq
 \Gamma' \int_{x_0}^t {\Phi^*}' \left( e^{\varepsilon(u-t)} \right) du
\leq
2 \Gamma' \Phi^*(1)\int_{x_0}^t  e^{\frac{\varepsilon}{1-\theta}  (u-t)}  du
\leq
\frac{2 (\theta-1) \Phi^*(1) \Gamma'}{\varepsilon}\cdot
$$
Hence for $t \in [x_0,x]$ the quantity $\eta \frac{f'_x(t)}{f_x(t)}$ is non-negative but bounded
away from zero.
So the value of $\Phi$ and its modification $H_\Phi$ on this quantity are comparable.
Consequently there exists a number $C$ (depending on $\Phi,\Gamma',\eta,\varepsilon,\theta$)
such that for $x\ge x_0$,
\begin{equation} \label{eq:alpha}
\widetilde \alpha_x^+
\leq
C \int_{x_0}^x  \Phi
\left( \eta \frac{f_x'}{f_x}  \right) \Phi'(f_x)^2 d\nu .
\end{equation}
At this stage, we need upper estimates for $f_x(t)$ and $c_x$.
 Integrating by parts as
in the proof of Lemma~\ref{lem:maj-norme} we get for $t\in (x_0,x)$
\begin{eqnarray*}
     f_x(t) &=& \int_{x_0}^t {\Phi^*}'\big(c_x e^{V(u)}\big)du \\
                     &=& \frac{\Phi^*\big(c_x e^{V(t)} \big)}{c_x e^{V(t)} V'(t)} -
                  \frac{\Phi^*\big(c_x e^{V(x_0)} \big)}{c_x e^{V(x_0)} V'(x_0)}+
                   \int_{x_0}^t \frac{ \Phi^*\big(c_x e^{V(u)}\big) }{c_x e^{V(u)}}
                   \left(1+\frac{V"(u)}{V'(u)^2} \right) \, du.
\end{eqnarray*}
Our choice of $x_0$ guarantees  $ -1\le \frac{V''(x)}{ V'(x)^2} \le 1/\theta$  for $x\ge x_0$.
Proceeding exactly as in the proof of Lemma~\ref{lem:maj-norme} yields
\begin{equation} \label{eq:fx}
f_x(t)  \leq F_x(t):=\theta(\theta-1) \frac{{\Phi^*}' \left( c_x e^{V_\nu(t)} \right)}{V'_\nu(t)} \cdot
\end{equation}
In order to estimate $c_x$, we use the above formula for $t=x$. Since $V$ is non-decreasing after $x_0$
and $\Phi^*(u)/u$ is also non-decreasing, we can write
\begin{eqnarray*}
    1= f_x(x)&  \ge&  \frac{\Phi^*\big(c_x e^{V(x)} \big)}{c_x e^{V(x)} V'(x)} -
                  \frac{\Phi^*\big(c_x e^{V(x_0)} \big)}{c_x e^{V(x_0)} V'(x_0)}+
                  \frac{\Phi^*\big(c_x e^{V(x_0)}\big) }{c_x e^{V(x_0)}}
                   \int_{x_0}^x
                   \left(1+\frac{V"(u)}{V'(u)^2} \right) \, du
    \\
              &=&   \frac{\Phi^*\big(c_x e^{V(x)} \big)}{c_x e^{V(x)} V'(x)} +
                  \frac{\Phi^*\big(c_x e^{V(x_0)} \big)}{c_x e^{V(x_0)}}
                    \left(x-\frac{1}{V'(x)} -x_0 +\frac{1}{V'(x_0)} -1\right).
\end{eqnarray*}
Recall that  $V'(x)\ge \varepsilon$ for $x\ge x_0$. Setting  $ x_1:=x_0+\varepsilon^{-1}+1$, we have obtained
 for $x\ge x_1$,
 $$ 1\ge  \frac{\Phi^*\big(c_x e^{V(x)} \big)}{c_x e^{V(x)} V'(x)} \ge
     \frac{{\Phi^*}'\big(c_x e^{V(x)} \big)}{2 V'(x)}=  \frac{{\Phi'}^{-1}\big(c_x e^{V(x)} \big)}{2 V'(x)},$$
hence
\begin{equation}
   \label{eq:cx}
   c_x \le e^{-V(x)} \Phi'\big(2 e^{V'(x)} \big).
\end{equation}
Now we go back to the estimate of $\widetilde{\alpha}_x^+$ given in \eqref{eq:alpha}. We give a
pointwize estimate
of the function in the integral of this equation: on $[x_0,x]$ it holds
\begin{eqnarray*}
  \Phi\left( \eta \frac{f'_x}{f_x}\right) \Phi'(f_x)^2
        &\le& \min(\eta^2,\eta^\theta)  \Phi\left( \frac{f'_x}{f_x}\right)\theta^2 \frac{\Phi(f_x)^2}{f_x^2}\\
        &\le& \min(\eta^2,\eta^\theta)\theta^2 \Gamma \Phi\left( f'_x\right) \frac{\Phi(f_x)}{f_x^2}\\
        &\le& \min(\eta^2,\eta^\theta)\theta^2 \Gamma f'_x\Phi'( f'_x) \frac{\Phi(F_x)}{F_x^2}\\
        &\le& \min(\eta^2,\eta^\theta)\theta^2 \Gamma f'_x\Phi'( f'_x) \frac{\Phi'(F_x)}{F_x},
\end{eqnarray*}
where we have used that $\Phi(x)/x^2$ is a non-decreasing function, together with the upper bound  $f_x\le F_x$ given in \eqref{eq:fx}.
In the following, $C_1,C_2,C_3$ are numbers depending on $\Phi, V,\eta,\theta,\varepsilon$ but not on $x$.
We also use repeatedly Lemma~\ref{lem:a2} to pull constants out of $\Phi$ or $\Phi'$. We get from \eqref{eq:alpha}, Lemma~\ref{lem:a3}
and the latter estimate
\begin{eqnarray*}
\widetilde \alpha_x^+
& \leq &
C_1 \int_ {x_0}^x  {\Phi'}^{-1}(c_x e^V)   c_x e^V \Phi'\left(\frac{{\Phi'}^{-1}(c_x e^V)}{V'} \right) \frac{V'}{{\Phi'}^{-1}(c_x e^V)} d\nu
 \\
& \leq & C_2 \int_ {x_0}^x     c_x e^V \frac{\Phi'\left({\Phi'}^{-1}(c_x e^V)\right)}{\Phi'(V')}  V' d\nu
\\
& =&
C_2 c_x^2
\int_{x_0}^x  \frac{V'(t) e^{V(t)}}{ \Phi' (V'(t))} dt .
\end{eqnarray*}
An integration by part formula leads to
\begin{eqnarray*}
\int_{x_0}^x  \frac{V'(t) e^{V(t)}}{ \Phi' (V'(t))} dt
& = &
\frac{e^{V(x)}}{\Phi' (V'(x))} -
\frac{e^{V(x_0)}}{\Phi' (V'(x_0))}
+
\int_{x_0}^x  \frac{V''(t)\Phi'' (V'(t))}
{{\Phi'}^2 (V'(t))}  e^{V(t)}dt \\
& \leq &
\frac{e^{V(x)}}{\Phi' (V'(x))}
+
\eta
\int_{x_0}^x  \frac{V''(t)}{{V'}^2(t)}
\frac{V'(t)e^{V(t)}}{\Phi' (V'(t))}  dt \\
& \leq &
\frac{e^{V(x)}}{\Phi' (V'(x))}
+
(1-\varepsilon)
\int_{x_0}^x \frac{V'(t)e^{V(t)}}{\Phi' (V'(t))}  dt
\end{eqnarray*}
where we have used
the assumption  $V'' / {V'}^2\le (1-\varepsilon)/\eta$ on $[x_0,+\infty)$. Hence for $x\ge x_0$
$$
\int_{x_0}^x  \frac{V'(t) e^{V(t)}}{ \Phi' (V'(t))} dt
\leq
\frac{1}{\varepsilon} \frac{e^{V(x)}}{\Phi' (V'(x))}  .
$$
Combining this bound with  the one on $\widetilde{\alpha}_x^+$ and the estimate \eqref{eq:cx} on $c_x$
gives, as claimed, that for $x\ge x_1$,
$$ \widetilde{\alpha}_x^+ \le C_3  e^{-V(x)} \Phi'(V'(x)).$$
\end{proof}
As an immediate consequence we get a converse statement to the criterion of Theorem~\ref{th:critphi-simple}

\begin{theorem}\label{th:critphi-simple-reciprocal}
Let $\Phi$ be a twice differentiable, strictly convex, nice Young function. Assume that on $\mathbb R^+$
the function $\Phi(x)/x^2$ is non-decreasing, the functions $\Phi(x)/x^\theta$ and $\Phi'(x)/x^\eta$
are non increasing for some $\theta>2,\eta>0$. Also assume that there exists $\Gamma\in\mathbb R$ such that
for all $x,y\ge 0$, $\Gamma \Phi(xy)\ge \Phi(x)\Phi(y)$.

Let $d\mu(x)=e^{-V(x)}dx$ be a probability measure on $\mathbb R$. Assume that $V$  of class $\mathcal C^2$
in neighborhoods of $+\infty $  such that:
\begin{enumerate}
\item $\displaystyle\liminf_{x\to +\infty } \mathrm{sign}(x) V'(x)>0$,
\item $\displaystyle -1<\liminf_{x\to +\infty } \frac{V''(x)}{V'(x)^2}  \le
                   \limsup_{x\to +\infty } \frac{V''(x)}{V'(x)^2} < \frac{1}{\max(\theta,\eta)}$,
\item there exists a constant $\kappa<+\infty $ such that for all locally Lipschitz $f$ on $\mathbb R$
$$\ent_{\mu} (f^2) \le \kappa \int_{\mathbb R}H_\Phi\left(\frac{f'}{f}\right)f^2 d\mu.$$
\end{enumerate}
Then   $\displaystyle\limsup_{x\to +\infty }\frac{V(x)}{\Phi\big(|V'(x)| \big)}<+\infty$ .

\end{theorem}
\begin{remark} A symmetric statement holds for $-\infty$.
\end{remark}
\begin{proof}
By Proposition~\ref{prop:min} for $x>m$ a median of $\mu$, it holds
  $$\mu\big([x,+\infty)\big) \log\left(1+\frac{1}{2 \mu\big([x,+\infty)\big) } \right) \le \kappa \widetilde\alpha_x^+.$$
  For $x$ large enough, Proposition~\ref{prop:al} provides an upper bound on $\widetilde\alpha_x^+$
  and $V$ is $\mathcal C^2$ so by Lemma~\ref{lem:equiv} the term $\mu\big([x,+\infty)\big)$ is lower bounded
  (and small enough to be where the function $t\log(1+1/(2t))$ increases). The conclusion follows easily.
 \end{proof}


\section{Concentration of measure phenomenon}

By Herbst argument, logarithmic Sobolev inequalities imply Gaussian
concentration, see e.g. \cite{Ane,ledoCMLS}.  Bobkov and Ledoux showed
that their modified inequality implies an improved form of exponential
concentration for products measures \cite{bobkl97pitc}, thus extending
a well-known result by Talagrand for the exponential measure
\cite{tala91niic}.  In this section we show that the argument may be
adapted to more general modified inequalities.

For a  convex function $H:[0,+\infty )\to \mathbb R^+$ we define
 $$\omega_H(x)=\sup_{t>0} \frac{H(tx)}{H(t)},\quad x\ge 0.$$
Clearly $\omega_H(0)=0$ and on $(0,+\infty)$ it is either identically infinite
or everywhere finite (exactly when $H$ satisfies the $\Delta_2$ condition).
One easily checks that $\omega_H\ge H/H(1)$ is convex and satisfies $\omega_H(ab)\le \omega_H(a)\omega_H(b)$
for all $a,b \ge 0$. Moreover if $H(x)/x^2$ is non decreasing for $x>0$ then so is the function
$\omega_H(x)/x^2$.

\begin{proposition} \label{prop:conc}
   Let $\mu$ be a probability measure on $\dR$ and $\mu^n$ the
   $n$-fold product measure on $\dR^n$.  Let $H:\mathbb R\to
   [0,+\infty ]$ be an even convex function. Assume that $x \mapsto H(x)/x^2$ is non-decreasing on $(0,+\infty )$.
   If
   there exists $\kappa<+\infty $ such that every locally Lipschitz $f : \dR
   \rightarrow \dR$ satisfies
\begin{equation} \label{eq:lsm}
\ent_{\mu}(f^2) \leq \kappa \int H \left( \frac{f'}{f} \right) f^2 d\mu,
\end{equation}
then every locally Lipschitz $F: \dR^n \rightarrow \dR$ with $\sum_{i=1}^n H(\partial_i
F) \leq a$ $\mu^n$-a.e. verifies
$$
\mu^n \left( \left\{F - \mu^n(F) \geq r \right\} \right) \leq e^{-K
  \omega_H^* \left( \frac{2r}{K} \right) } \qquad \forall r \geq 0
$$
where $\omega_H^*$ is the conjugate of $\omega_H$ and $K=a\kappa$.
\end{proposition}

\begin{proof}
 We may assume that $\omega_H$ is everywhere finite otherwise there is nothing to prove.
   Fix $F : \dR^n \rightarrow \dR$ with $\sum_{i=1}^n H(\partial_i F)
   \leq a$. Assume first that $F$ is integrable. By tensorisation of the modified logarithmic Sobolev
   Inequality \eqref{eq:lsm} (see \cite{gentilgm}), any locally Lipschitz $f:
   \dR^n \rightarrow \dR$ verifies
   $$
   \ent_{\mu^n}(f^2) \leq \kappa \int \sum_{i=1}^n H \left(
      \frac{\partial_i f}{f} \right) f^2 d\mu^n .
   $$
   Plugging $f:=e^{\frac{\lambda}{2}F}$, $\lambda \in \dR^+$, leads
   to
\begin{eqnarray*}
\ent_{\mu^n}(e^{\lambda F})
& \leq &
\kappa \int \sum_{i=1}^n H \left( \frac{\lambda}{2}\partial_i F \right) e^{\lambda F} d\mu^n \\
& \leq  &
 \kappa\,a\, \omega_H \left( \frac{\lambda}{2} \right) \int e^{\lambda F} d\mu^n  .
\end{eqnarray*}
Define $\Psi(\lambda):=\int e^{\lambda F} d\mu^n$. Then
$\ent_{\mu^n}(e^{\lambda F})=\lambda \Psi'(\lambda) - \Psi(\lambda)
\log\Psi(\lambda)$.  Hence, by definition of $K$,
$$
\lambda \Psi'(\lambda) - \Psi(\lambda) \log\Psi(\lambda) \leq K \omega_H
\left( \frac{\lambda}{2} \right) \Psi( \lambda) \qquad \forall \lambda
\geq 0 .
$$
In particular, dividing by $\lambda^2 \Psi(\lambda)$,
$$
\frac{d}{d\lambda} \left( \frac{\log \Psi(\lambda)}{\lambda}
\right) \leq K \frac{\omega_H \left( \frac{\lambda}{2} \right)}{\lambda^2}
\qquad \forall \lambda >0 .
$$
Note that $\lim_0 \frac{\log \Psi(\lambda)}{\lambda} = \mu^n(F)$.
Hence integrating leads to
$$
\int e^{\lambda (F-\mu^n(F))} d\mu^n \leq \exp\left\{ K \lambda
   \int_0^\lambda \frac{\omega_H \left( \frac{u}{2} \right)}{u^2} du \right\}
.
$$
Chebichev Inequality finally gives for any $r\geq0$, any $\lambda
>0$,
\begin{eqnarray*}
\mu^n\left( \left\{ F - \mu^n(F) \geq r \right\} \right)
& \leq &
e^{-\lambda r} \int e^{\lambda(F-\mu^n(F))} d\mu^n
\end{eqnarray*}
which leads to
$$
\mu^n\left( \left\{ F - \mu^n(F) \geq r \right\} \right) \leq \exp
\left\{ - K \sup_{\lambda >0} \left[ \frac{2r}{K} \frac{\lambda}{2}
      -\lambda \int_0^\lambda \frac{\omega_H \left( \frac{u}{2} \right)}{u^2}
      du \right] \right\} .
$$
The conclusion follows from the inequality
$$
\lambda \int_0^\lambda \frac{\omega_H \left( \frac{u}{2} \right)}{u^2} du
\leq \omega_H \left( \frac{\lambda}{2} \right),
$$
which is proved as follows: let $\theta (\lambda):= \int_0^\lambda
\frac{\omega_H \left( u/2 \right)}{u^2} du $.  The result is equivalent to
$\theta(\lambda) \leq \lambda \theta'(\lambda)$.  Now since
$\theta'(\lambda)=\frac{\omega_H \left( \lambda / 2 \right)}{\lambda^2} =
\frac 14 \frac{\omega_H \left( \lambda / 2
  \right)}{\left(\lambda/2\right)^2}$ is non decreasing, $\theta$ is
convex.  In turn, since $\theta(0)=0$, $\theta(\lambda) \leq \lambda
\theta'(\lambda)$ as expected.

The proof is complete for $F$ integrable. A standard truncation argument, see e.g. \cite[Lemma 7.3.3]{Ane},
shows that $F$ is automatically integrable.
\end{proof}

\begin{theorem} \label{th:conc}
   Let $\mu$ be a probability measure on $\dR$, which we assume to be absolutely continuous with
  respect to Lebesgue's measure. Let  $H:\mathbb R\to \mathbb R^+$ be an even
    convex function, with $H(0)=0$.
    Assume  that $x \mapsto H(x)/x^2$ is non-decreasing  for $x>0$ and that $H^*$ is strictly convex. If
    there exists $\kappa<+\infty $ such that every
   locally Lipschitz $f : \dR \rightarrow \dR$ satisfies
\begin{equation*}
\ent_{\mu}(f^2) \leq \kappa \int H \left( \frac{f'}{f} \right) f^2 d\mu,
\end{equation*}
then every Borel set $A \subset \dR^n$ with $\mu^n(A) \geq \frac 12$
satisfies
$$
1 - \mu^n \left( A + \Big\{x:\sum_{i=1}^n H^*(x_i) < r
   \Big\} \right) \leq e^{-Kr} \qquad \forall r \geq 0
$$
where $K=\omega_ H(2)\, \kappa \, \omega_H^* \left( \frac{1}{\omega_H(2)\, \kappa}
\right)$.
\end{theorem}
\begin{remark}
   The hypothesis of strict convexity of $H^*$ is here for technical reasons. In practice
   $H^*$ often fails to be strictly convex on a set $[a,b]\subset (0,+\infty )$. In this case
it is easy to build an even strictly convex function $I\ge H^*$ which actually coincides with
$H^*$ outside of a slightly larger interval and satisfies $I_r' \le 2 {H^*}_r'$ on $\mathbb R^+$.
Following the proof of the theorem with $I$ instead of $H^*$ then yields the concentration inequality
claimed in the above theorem, only with a worse constant.
\end{remark}

\begin{proof}
   We start with establishing a useful inequality verified by $H$. Since $H(x)/x^2$ is non-decreasing on $(0,+\infty)$
   it follows that $H^*(x)/x^2$ is non-increasing on this interval, and taking right derivatives that
    $2H^*(x) \ge x(H^*)'_r(x)$ for $x>0$ (actually Lemma~\ref{lem:a1} is valid without differentiability).
    Next we use the easy inequality $H^*(x)\ge H\big(H^*(x)/x \big)$ for $x>0$ (it is usually written
    in the following nicer but more restrictive form $H^{-1}(x){H^*}^{-1}(x)\ge x$). It follows that
    \begin{equation}\label{eq:spurious}
    H^*(x) \ge H\left(\frac{(H^*)'_r(x)}{2} \right) \ge \frac{1}{\omega_H(2)} H\big((H^*)'_r(x)\big).
    \end{equation}
   Let $A \subset \dR^n$ with $\mu^n(A) \geq \frac 12$ and $F_A(x)=
   \inf_{z \in A} \sum_{i=1}^n H^*(x_i -z_i)$ for
   $x=(x_1,\cdots,x_n) \in \dR^n$. For $r >0$ set further
   $F=\min(F_A,r)$.
   We claim that Lebesgue a.e and thus $\mu^n$-a.s., it holds
 \begin{equation}\label{eq:H}
   \sum_{i=1}^n H(\partial_i F) \leq \omega_H(2) \,r.
  \end{equation}
 First let us develop the consequence of this claim.
   Note that $F_A=0$ on $A$. Thus, $\int F d\mu^n \leq
   r(1-\mu^n(A)) \leq \frac r2$.  Hence, since $\left\{ F \geq r
   \right\} \subset \left\{ F-\mu^n(F) \geq \frac r2 \right\}$,
   Proposition \ref{prop:conc} ensures that
$$
\mu^n \left( \left\{ F \geq r \right\} \right)
 \leq
\mu^n \left( \left\{ F -\mu^n(F) \geq \frac r2  \right\}  \right)
\leq
\exp \left\{ - \omega_H(2)\, r\,\kappa\,  \omega_ H^* \left( \frac{1}{\omega_H(2) \kappa} \right) \right\} .$$
This leads to the expected result since one can easily see that
$$\left\{ F < r \right\} = \left\{ F_A < r \right\} \subset A +
\Big\{x:\sum_{i=1}^n H^*(x_i) < r \Big\}\cdot$$
Finally we establish the claim \eqref{eq:H}. Since $H^*$ is convex and always finite, it is locally Lipschitz
and one easily checks that this property passes to $F$. Hence $F$ is almost everywhere differentiable and
the set $\{x;\; \nabla F(x)\neq 0 \;\mbox{and}\; F=r\}$ is negligeable. Hence we may restrict to points where $F<r$
and thus $F=F_A<r$ and $F_A$ is differentiable. Denote $\mathcal H (x)=\sum_{i=1}^{n} H^*(x_i)$.

 We shall first
prove that when $F_A$ is differentiable at $x$, there exits a unique $a\in \overline{A}$ such that
$F_A(x)=\mathcal H(x-a)$. Assume that $F_A$ is differentiable at $x$ and that
there exist $a\neq b$ in $\overline{A}$ such that $F_A(x)=\min_{c\in \overline{A}}\mathcal H(x-c)
=\mathcal H(x-a)=\mathcal H(x-b)$.
Consider the function
$L:[0,1]\to\mathbb R$ defined by $L(u)= \mathcal H\big(x-(ua+(1-u)b)\big)$. Since it is strictly convex and $L(0)=F_A(x)=L(1)$
it follows that $L'_r(0)<0<L'_\ell(1)$. Since $b\in \overline{A}$ it holds for $t\in [0,1]$
$$ F_A(x+t(b-a))\le \mathcal H(x+t(b-a)-b)=L(t),$$
with equality at $t=0$.
It follows that  $DF_A(x).(b-a) \le L'_r(0)<0$.
On the other hand, since $a\in \overline{A}$, it holds for $t\in[-1,0]$,
 $$ F_A(x+t(b-a))\le \mathcal H(x+t(b-a)-a)=L(1+t),$$
with equality at $t=0$.
It follows that  $DF_A(x).(b-a) \ge L'_\ell(1)>0$ which contradicts our previous bound.

To complete the proof of the claim, we consider a point $x$ where $F_A$ is differentiable and $F_A(x)<r$
and we consider $a\in \overline{A}$ the unique minimizer for $\mathcal H(x-\cdot)$ on $\overline{A}$.
An easy consequence of the uniqueness is that for every sequence $y^k$ converging to $x$
and $a^k\in \overline{A}$ such that $F_A(y^k)= \mathcal H(y^k-a^k)$, the sequence $a^k$ converges to $a$.
Let $t_k$ be a sequence of positive numbers converging to zero.
Then, denoting by $e^i$ the $i$-th vector in the canonical basis of $\mathbb R^n$,
\begin{eqnarray*}
   F_A(x+t_k e^i)- F_A(x) &=& \inf_{c\in \overline A} \mathcal H (x+t_ke^i-c) -\mathcal H(x-a) \\
    &\le & \mathcal H (x+t_ke^i-a) -\mathcal H(x-a) =H^*(x_i +t_k-a_i)-H^*(x_i-a_i).
\end{eqnarray*}
Dividing by $t_k>0$ and taking limits yields $\partial_i F_A(x) \le
{H^*_r}'(x_i-a_i) \le {H^*_r}'(|x_i-a_i|).$ Similarly, if we denote
by $a^k$ a minimizer of $c\in \overline A \mapsto \mathcal H(x+t_k
e^i-c)$
\begin{eqnarray*}
    F_A(x+t_k e^i)- F_A(x) &=&  \mathcal H (x+t_ke^i-a^k) -\inf_{c\in \overline A}\mathcal H(x-c) \\
    &\ge & \mathcal H (x+t_ke^i-a^k) - \mathcal H (x-a^k) = H^*(x_i +t_k-a_i^k)-H^*(x_i-a_i^k )\\
    &\ge & t_k {H_r^*}'(x_i-a^k_i),
\end{eqnarray*}
by convexity. Recall that $a^k$ converges to $a$. Hence letting $k$ to infinity we get $\partial_i F_A(x) \ge
{H_\ell^*}'(x_i-a_i)\ge -{H_r^*}'(|x_i-a_i|)$.
Eventually when $F_A(x)=F(x)<r$
\begin{eqnarray*}
\sum_{i=1}^{n} H\big(\partial_i F(x) \big) &\le&  \sum_{i=1}^{n} H\big({H_r^*}'(|x_i-a_i|)  \big)\le
  \omega_H(2) \sum_{i=1}^{n} H^*(x_i-a_i)\\
 &=& \omega_H(2) \mathcal H(x-a) =\omega_H(2) F_A(x) <\omega_H(2) r,
\end{eqnarray*}
using \eqref{eq:spurious} and the definition of $a$ as a minimizer.
\end{proof}

If $H=H_\Phi$ is the modification of an even convex $\Phi:\mathbb
R\to \mathbb R^+$ with $\Phi(x)/x^2$ non-decreasing on $\mathbb R^+$
one easily checks that there exists $x_0$ such that $H^*_\Phi(x)$ is
comparable to $x^2$ up to multiplicative constants if $|x|\le x_0$,
and $H_\Phi^*(x)=\Phi^*(x)$ otherwise. Then, separating coordinates
$x_i$ of absolute value less or more than $x_0$, one gets that
 there exists a constant $c$
(depending on $\Phi$) such that for any $r$,
$$
\left\{x:\sum_{i=1}^n H_\Phi^*(x_i) < r \right\}
\subset
\sqrt{cr} B_2 + \left\{x : \sum_{i=1}^n \Phi^*(x_i) < cr \right\}.
$$
Let $\omega_{\Phi^*}(t):=\sup_{x >0} \frac{\Phi^*(tx)}{\Phi^*(x)}$ for
$t >0$ and $B_{\Phi^*} :=\left\{x: \sum_{i=1}^n \Phi^*(x_i) < 1
\right\}$.  For any $x$ such that $\sum_{i=1}^n \Phi^*(x_i) < s$, we
have
$$
\sum_{i=1}^n \Phi^*\left( \omega_{\Phi^*}^{-1}\left (\frac
      1s\right) x_i \right) \leq \omega_{\Phi^*} \left(
   \omega_{\Phi^*}^{-1}\left (\frac 1s \right) \right) \sum_{i=1}^n
\Phi^*(x_i) < 1 .
$$
Thus $\left\{x:\sum_{i=1}^n \Phi^*(x_i) < s \right\} \subset
\frac{1}{\omega_{\Phi^*}^{-1}\left (\frac 1s\right)} B_{\Phi^*}$.
Hence, under the hypotheses  of Theorem \ref{th:conc} we have
for any Borel set $A \subset \dR^n$ with $\mu^n(A) \geq \frac 12$,
\begin{equation}\label{eq:conc}
\mu^n \left( A +
\sqrt{r} B_2 + \frac{1}{\omega_{\Phi^*}^{-1}\left (\frac{1}{r} \right)} B_{\Phi^*}
\right) \ge \mu^n\left(A+  \Big\{x:\sum_{i=1}^n H_\Phi^*(x_i) < r \Big\}\right)
\geq 1-e^{-Cr} \qquad \forall r \geq 0
\end{equation}
for some constant $C$ independent on $r$.
Such concentration inequalities were established by Talagrand  \cite{tala91niic,tala95cmii} for the exponential measure and
later for even log-concave measures,
via inf-convolution inequalities (which are strongly related to transportation cost inequalities).
More recently Gozlan derived such inequalities from his criterion for transportation inequalities
on the line \cite{gozl06ctti}.
We conclude this section with concrete examples.


\begin{example}
   Let $\Phi_q(x)=|x|^q$,
   $q \geq 2$ and $H_q(x)=H_{\Phi_q}(x)=\max(x^2,|x|^q)$.  Straightforward calculations give
   $$
   H_q^*(x) = \left\{ \begin{array}{ll}
         x^2/4 & \mbox{if } x \leq 2 \\
         x-1 & \mbox{if } 2 \leq x \leq q \\
         (q -1) \left( x/q \right)^\frac{q}{q -1}
         & \mbox{if } x \geq q
\end{array}
\right. .
$$
Here $\omega_{\Phi_q^*}=\Phi_q^* =
C_q |x|^{q^*}$ with $\frac 1 q + \frac{1}{q^*} =1$.
Let $B_{q^*} := \left\{ x: \sum_{i=1}^n |x_i|^{q^*} <1
\right\}$ be the $\ell^{q^*}$-unit ball in $\dR^n$.  If $\mu$
satisfies the modified logarithmic Sobolev Inequality
\eqref{eq:phialpha}, there exists a constant $C'_q$ (depending
only on $q$) such that
$$
1 - \mu^n \left( A + \sqrt{r} B_2 + r^\frac{1}{q^*} B_{q^*}
\right) \leq e^{-C'_q r} \qquad \forall r \geq 0
$$
for any $A$ with $\mu^n(A) \geq \frac 12$.  In particular, thanks
to Corollary \ref{cor:mualpha}, the measures $d\mu_\beta(x)=
Z_\beta^{-1} e^{-|x|^{\beta}} dx$ satisfy the latter concentration
result for any $\beta \geq q^*>1$.

Note that the limit case $q^*=1$ or $q=+\infty$ is not treated in our argument.
It corresponds to the case when $H(x)=x^2 \ind_{|x|<c}+\infty  \ind_{|x|\ge c}$
treated by Bobkov and Ledoux \cite{bobkl97pitc}. Our ``extension'' does not cover this
case since for technical reasons we considered only functions $H$ taking finite values.
On the other hand combining Corollary~\ref{cor:gen} with the above theorem and remark,
yields similar concentration properties for a wide class of even log-concave measures
with an intermediate behaviour between exponential and Gaussian.
\end{example}

\section{Appendix on Young functions} \label{sec:orlicz}

In this section we collect some useful results and definition on Orlicz spaces.
We refer the reader to \cite{rao-ren} for demonstrations and complements.
\begin{definition}[Young function]
A function $\Phi : \dR \rightarrow [0,\infty]$ is a \emph{Young function}
if it is convex, even, such that $\Phi(0)=0$, and
$\lim_{x \rightarrow +\infty} \Phi(x)=+\infty$.
\end{definition}
The Legendre transform $\Phi^*$ of $\Phi$ is defined by
$\Phi^*(y)= \sup_{x \geq 0} \{ x|y| - \Phi(x) \}$. It is
a lower semi-continuous Young function called the \emph{complementary function} or \emph{conjugate} of $\Phi$.
Among the Young functions, we call  \emph{nice Young function} those
which take only  finite values and such that $\Phi(x)/x \rightarrow \infty$
as $x \rightarrow \infty$,
 $\Phi(x)=0 \Leftrightarrow x=0$ and
$\Phi'(0)=0$.

For any nice Young function $\Phi$, the conjugate of $\Phi^*$ is $\Phi$
and for any $x>0$,
$$
x \leq \Phi^{-1}(x) (\Phi^*)^{-1}(x) \leq 2x .
$$
The simplest example of nice Young function is $\Phi(x)=\frac{|x|^p}{p}$,
$p>1$, for which, $\Phi^*(x)=\frac{|x|^q}{q}$, with $1/p+1/q=1$.

Now let $(\cal X,\mu)$ be a measurable space, and $\Phi$ a Young function.
The space
$$
\mathbb{L}_\Phi(\mu)=\{f: {\cal X} \rightarrow \dR \mbox{ measurable};
\exists \alpha > 0,
\int_{\cal X} \Phi(\alpha f) < + \infty \}
$$
is called the {\it Orlicz space} associated to $\Phi$.
When $\Phi(x)=|x|^p$, then $\mathbb{L}_\Phi(\mu)=\mathbb{L}^{p}(\mu)$, the standard
Lebesgue space.
There are two natural equivalent norms which give to $\mathbb{L}_\Phi(\mu)$ a structure of
Banach space. Namely
$$
\| f \|_{\Phi} = \inf\{ \lambda >0; \int_{\cal X}
\Phi \left( \frac{f}{\lambda} \right) d\mu \leq 1 \}
$$
and
$$
N_\Phi(f)=\sup \{\int_{\cal X} |fg|d\mu; \int_{\cal X}
\Phi^*(g) d\mu \leq 1 \} \;.
$$
Note that we invert the notation with respect to \cite{rao-ren}.
For $\alpha\in [1,\infty]$ we denote the dual coefficient $\alpha^*\in [1,\infty ]$. It is defined
by the equality $\frac{1}{\alpha}+\frac{1}{\alpha^*}=1$.

\begin{lemma}
   \label{lem:a1} Let $\alpha\in (1,+\infty )$.
 Let $\Phi$ be a differentiable, strictly convex nice Young function.
  Then the following assertions are equivalent:
  \begin{enumerate}
  \item The function $\Phi(x)/x^\alpha$ is non-decreasing for $x>0$.
  \item For $x\ge 0$, $x\Phi'(x)\ge \alpha \Phi(x).$
  \item For $x\ge 0$, $x{\Phi^*}'(x)\le \alpha^* \Phi^*(x).$
    \item The function $\Phi^*(x)/x^{\alpha^*}$ is non-increasing for $x>0$.
  \end{enumerate}
\end{lemma}
Note that $\Phi$ and $\Phi^*$ play symmetric roles so that similar equivalent formulations exist for the
property: $\Phi(x)/x^\alpha$  is non-increasing for $x\ge 0$.
\begin{proof}
Plainly, the first two statements are equivalent by taking derivatives, and the last two as well.
We show that $(ii)$ implies $(iii)$. 
Our hypotheses ensure that $\Phi'$ is a bijection of $[0;+\infty)$; its inverse is ${\Phi^*}'$.
Since  for $x\ge 0$, $x\Phi'(x)\ge \alpha \Phi(x)$,
\begin{eqnarray*}
  \Phi^*(x) &=& \sup_y\big\{ xy-\Phi(y) \big\}=x {\Phi'}^{-1}(x) -\Phi\big( {\Phi'}^{-1}(x)\big) \\
      & \ge & x {\Phi'}^{-1}(x) - \frac{1}{\alpha} {\Phi'}^{-1}(x) \Phi'\big( {\Phi'}^{-1}(x)\big)\\
      &=& \left(1- \frac{1}{\alpha} \right) x {\Phi'}^{-1}(x).
\end{eqnarray*}
Hence using that $\Phi'$ and ${\Phi^*}'$ are inverse function,
\begin{equation}\label{eq:phitheta}
\Phi^*(x) \ge \left(1- \frac{1}{\alpha} \right) x {\Phi^*}'(x)=\frac{1}{\alpha^*}  x {\Phi^*}'(x).
\end{equation}
A similar argument yields the converse implication.
\end{proof}

The next lemma is obvious, but convenient.
\begin{lemma}
   \label{lem:a2} Let $0<\alpha<\theta$.
 Let $\Phi$ be a  differentiable function on $[0,+\infty )$ such that the function
  $\Phi(x)/x^\alpha$ is non-decreasing and $\Phi(x)/x^\theta$ is non-increasing.
 Then for $x > 0$, and $t\ge 1$ it holds
   $$  \Phi(tx)\le t^\theta \Phi(x),
    \quad \Phi'(tx)\le \theta t^{\theta-1} \frac{\Phi(x)}{x} \le \frac{\theta}{\alpha} t^{\theta-1}\Phi'(x).$$
 For for $x > 0$, and $t\in (0,1]$ it holds
   $$  \Phi(tx)\le t^\alpha \Phi(x),
    \quad \Phi'(tx)\le \theta t^{\alpha-1} \frac{\Phi(x)}{x} \le \frac{\theta}{\alpha} t^{\alpha-1}\Phi'(x).$$

\end{lemma}

\begin{lemma}
   \label{lem:a3} Let $1<\alpha< \theta$.
 Let $\Phi$ be a strictly convex differentiable nice Young function such that $\Phi(x)/ x^\alpha$ is non-decreasing
  for $x>0$ and $\Phi(x)/x^\theta$ is non-increasing for $x>0$. Assume that there exists
  $\Gamma \in \mathbb R^+$ such that for all $x,y\ge 0$ it holds
   $$ \Gamma \Phi(xy)\ge \Phi(x)\Phi(y).$$
  Then there exist real numbers $\Gamma_1,\Gamma_2,\Gamma_3\in \mathbb R^+$ such that
 for all $x,y\ge 0$,
   $$  \Gamma_1 \Phi'(xy)\ge \Phi'(x)\Phi'(y),
   \quad  {\Phi^*}'(xy)\le  \Gamma_2 {\Phi^*}'(x){\Phi^*}'(y),
   \quad  \Phi^*(xy)\le \Gamma_3\Phi^*(x)\Phi^*(y).  $$
\end{lemma}
\begin{proof}
   It is enough to deal with $x,y>0$. Our assumption and Lemma~\ref{lem:a1} allow to write
   $$  \frac{ \Gamma}{\alpha} \Phi'(xy)\ge \Gamma\frac{\Phi(xy)}{xy}
 \ge\frac{\Phi(x)}{x}\frac{\Phi(y)}{y} \ge \frac{1}{\theta^2} \Phi'(x) \Phi'(y),$$
   which gives the result for $\Phi'$ with $\Gamma_1=\theta^2\Gamma/\alpha$.
  Applying the inequality for $\Phi'$ to $x={\Phi^*}'(a),y={\Phi^*}'(b)$ and since ${\Phi^*}'$ is the
   inverse bijection of $\Phi'$ we get
    $$ {\Phi^*}'(a){\Phi^*}'(b)\ge {\Phi^*}'\Big(\frac{1}{\Gamma_1} ab\Big).$$
 Combining the hypotheses on the growth of $\Phi$ with Lemma~\ref{lem:a1} and Lemma~\ref{lem:a2} we obtain
  that for all $x,t>0$,
   $$ {\Phi^*}'(tx) \le \frac{\alpha^*}{\theta^*} \max\Big(t^{\alpha^*-1}, t^{\theta^*-1} \Big){\Phi^*}'(x) .$$
  Applying this inequality to $x= ab/\Gamma_1$ and $t=\Gamma_1$ shows that there exists $\Gamma_2>0$
  such that ${\Phi^*}'(ab/\Gamma_1)\ge {\Phi^*}'(ab)/\Gamma_2$. Hence the claimed inequality is valid for ${\Phi^*}'$.
  Finally
  $$  (\alpha^*)^2 \Gamma_2\frac{\Phi^*(x)\Phi^*(b)}{ab} \ge
  \Gamma_2 {\Phi^*}'(a){\Phi^*}'(b)\ge {\Phi^*}'(ab)\ge \theta^* \frac{\Phi^*(ab)}{ab},$$
  and the proof is complete.
\end{proof}




\noindent
{\it Mathematics Subject Classification}: 26D10, 60E15.\\
{\it Keywords}: Sobolev inequalities, concentration.
\bigskip

\noindent
F. B.: Institut de Mathématiques, Université Paul Sabatier, 31062 Toulouse cedex 09, FRANCE.
E-mail: barthe@math.ups-tlse.fr

\medskip
\noindent
C. R.: Laboratoire d'Analyse et Mathématiques Appliquées- UMR 8050,
 Universités de Marne la Vallée et de Paris 12-Val-de-Marne,
 Boulevard Descartes, Cité Descartes, Champs sur Marne,
 77454 Marne la Vallée Cedex 2, FRANCE.
E-mail: cyril.roberto@univ-mlv.fr

 \end{document}